\RequirePackage{ifpdf}
\ifpdf 
\documentclass[pdftex]{sigma}
\else
\documentclass{sigma}
\fi

\numberwithin{equation}{section}

\newtheorem{Theorem}{Theorem}[section]
\newtheorem{Corollary}[Theorem]{Corollary}
\newtheorem{Lemma}[Theorem]{Lemma}
\newtheorem{Proposition}[Theorem]{Proposition}
{\theoremstyle{definition}
\newtheorem{Note}[Theorem]{Note}
\newtheorem{Definition}[Theorem]{Definition}
\newtheorem{Problem}[Theorem]{Problem}
}

\begin{document}

\allowdisplaybreaks

\renewcommand{\PaperNumber}{099}

\FirstPageHeading

\ShortArticleName{The Universal Askey--Wilson Algebra}

\ArticleName{The Universal Askey--Wilson Algebra\\
and
the Equitable Presentation of
$\boldsymbol{U_q(\mathfrak{sl}_2)}$}

\Author{Paul TERWILLIGER}

\AuthorNameForHeading{P.~Terwilliger}

\Address{Department of Mathematics, University of Wisconsin, Madison, WI 53706-1388, USA}

\Email{\href{mailto:terwilli@math.wisc.edu}{terwilli@math.wisc.edu}}

\ArticleDates{Received July 19, 2011, in f\/inal form October 10, 2011;  Published online October 25, 2011}

\Abstract{Let $\mathbb F$ denote a f\/ield, and f\/ix a nonzero $q \in \mathbb F$ such
that $q^4\not=1$.
The universal Askey--Wilson algebra is the
associative $\mathbb F$-algebra $\Delta=\Delta_q$ def\/ined by generators and
relations in the following way.
The generators are
$A$, $B$, $C$. The relations assert that each~of
\[
A + \frac{qBC-q^{-1}CB}{q^2-q^{-2}},
\qquad
B + \frac{qCA-q^{-1}AC}{q^2-q^{-2}},
\qquad
C + \frac{qAB-q^{-1}BA}{q^2-q^{-2}}
\]
is central in $\Delta$.
In this paper we discuss a connection between
$\Delta$ and the $\mathbb F$-algebra
$U=U_q(\mathfrak{sl}_2)$.
To summarize the connection,
let $a$, $b$, $c$ denote mutually commuting indeterminates and
let
$\mathbb F \lbrack a^{\pm 1}, b^{\pm 1}, c^{\pm 1}\rbrack$
denote the $\mathbb F$-algebra
of Laurent polynomials in
$a$, $b$, $c$ that have all coef\/f\/icients in
$\mathbb F$.
We display
an injection of $\mathbb F$-algebras
$ \Delta \to
U \otimes_\mathbb F \mathbb F \lbrack a^{\pm 1}, b^{\pm 1}, c^{\pm 1}\rbrack
$.
For this injection we give the image of
$A$, $B$, $C$
and the
above three central elements,
in terms of the
equitable generators for~$U$.
The algebra~$\Delta $ has another central element of interest,
called the Casimir element~$\Omega$.
One signif\/icance of~$\Omega$ is the following.
It is known that
the center
of $\Delta$ is generated by
$\Omega$
 and the above three central
 elements,
 provided that  $q$ is not a root of unity.
 For the above injection
 we give the  image of
 $\Omega$
 in terms of the equitable generators for
 $U$.
 We also use the injection to
 show that
 $\Delta$ contains no zero divisors.}


\Classification{33D80; 33D45}

\vspace{-3mm}

\section{Introduction}

The Askey--Wilson polynomials
were introduced in \cite{awpoly} and soon became renown for
their algebraic, analytic, and combinatorial properties
\cite{ismail,koeswa}.
In his study
\cite{zhidd}
of the ``hidden symmetry'' of
these
polynomials,
A. Zhedanov introduced the Askey--Wilson
algebra AW(3).
This algebra is def\/ined
by generators and relations. The relations involve a nonzero
parameter $q$ and 5 additional parameters.
The algebra
is
 inf\/inite dimensional
and noncommutative.
Zhedanov's original presentation
involves three generators;  however one generator is
a $q$-commutator of the other two and is sometimes
eliminated. The remaining two generators
satisfy a pair of relations often called
the Askey--Wilson relations \cite{aw}.
These relations describe the
Askey--Wilson polynomials in
the following way.
Let $\lbrace p_n\rbrace_{n=0}^\infty$
denote a sequence of
Askey--Wilson polynomials
in a~variable~$\lambda$.
These polynomials
 are eigenvectors for
a certain $q$-dif\/ference operator, known
as the Askey--Wilson operator and denoted for the moment by $A$.
Let $B$ denote the linear operator
that sends $p(\lambda) \mapsto \lambda p(\lambda)$
for all polynomials $p(\lambda)$.
The operator $B$ acts on the basis
 $\lbrace p_n\rbrace_{n=0}^\infty$
 in an irreducible tridiagonal fashion,
 ref\/lecting the fact that
 $\lbrace p_n\rbrace_{n=0}^\infty$ satisfy
 a
 three-term recurrence.
In \cite{zhidd} Zhedanov
showed that $A$, $B$ satisfy a pair of Askey--Wilson relations.

Although the Askey--Wilson relations
are slightly complicated,
over time their signif\/icance became clear
as they found applications to
integrable systems
\cite{bas2,lavr,wz},
quantum groups
  \cite{GYZlinear},
linear algebra
\cite{LS99,aw},
quantum mechanics
\cite{odake},
and the double af\/f\/ine Hecke algebra
\cite{koo1,
dahater}.
We now describe the two
applications most  relevant to the present paper.

Our f\/irst application concerns a linear algebraic object
called a Leonard pair. This is a
 pair of diagonalizable linear transformations
on a f\/inite-dimensional vector space, each of which
acts in an irreducible tridiagonal fashion on an eigenbasis for
the other one
\cite[Def\/inition 1.1]{LS99}.
In~\cite{aw} Vidunas and the present author showed that
a Leonard pair satisf\/ies a pair of
Askey--Wilson relations.
This example is closely related to the
one involving the Askey--Wilson
polynomials.  By
\cite[Section~5]{TLT:array}
the Leonard pairs correspond
to a family
of orthogonal polynomials consisting of the $q$-Racah
polynomials and their relatives. The polynomials in this
family
are
special or limiting cases of the Askey--Wilson
polynomials~\cite{koeswa}.

Our second application is about quantum groups.
Consider the algebra
$U_q(\mathfrak{sl}_2)$
with the usual Chevalley generators
$e$, $f$, $k^{\pm 1}$~\cite{jantzen}.
Let $A$ denote an arbitrary
linear combination of $ek^{-1}$,~$f$,~$k^{-1}$ and let $B$ denote
 an arbitrary linear
combination of $e$, $fk$,  $k$. Then according to
   Granovski{\u\i} and Zhedanov~\cite{GYZlinear}
the elements $A$, $B$ satisfy
a pair of Askey--Wilson relations.
In~\cite{wz} Wiegmann and Zabrodin
extended this result by
displaying an element
$C$ in
$U_q(\mathfrak{sl}_2)$
such that
\begin{gather}
qAB-q^{-1}BA  =  g_C C + h_C,
\label{eq:zsym1}
\\
qBC-q^{-1}CB  =  g_A A + h_A,
\label{eq:zsym2}
\\
qCA-q^{-1}AC  =  g_B B + h_B,
\label{eq:zsym3}
\end{gather}
where $g_A$, $g_B$, $g_C$ and $h_A$, $h_B$, $h_C$
are scalars in the underlying f\/ield.
The equations
(\ref{eq:zsym1})--(\ref{eq:zsym3})
are often called the
$\mathbb Z_3$-symmetric Askey--Wilson relations
\cite{dahater}.
Upon eliminating $C$ in~(\ref{eq:zsym2}),~(\ref{eq:zsym3}) using
(\ref{eq:zsym1}) we obtain the Askey--Wilson relations
in the variables~$A$,~$B$. Upon substituting
$C'=g_C C + h_C$ in
(\ref{eq:zsym1})--(\ref{eq:zsym3}) we recover the
original presentation for AW(3) in the variables~$A$,~$B$,~$C'$.

We now recall the universal
Askey--Wilson algebra~$\Delta$~\cite{uaw}.
To motivate this algebra
consider the
relations
(\ref{eq:zsym1})--(\ref{eq:zsym3}). They
are attractive but one might object that
there are too many parameters.
To accomodate this objection
we will
eliminate
all the parameters besides~$q$.
We will do this
without
signif\/icantly reducing
the generality of the algebra
(although
we allow a minor technical assumption).
We f\/irst eliminate
 $g_A$, $g_B$, $g_C$ with the following
 change of variables.
Assume that
each of $g_Ag_B$, $g_Bg_C$, $g_Cg_A$ is a nonzero
square in the underlying f\/ield, and that $q^4\not=1$.
Now in
(\ref{eq:zsym1})--(\ref{eq:zsym3})
 replace $A$, $B$, $C$ by
\begin{gather*}
\frac{A f_A}{q^{-2}-q^2},
\qquad
\frac{B f_B}{q^{-2}-q^2},
\qquad
\frac{C f_C}{q^{-2}-q^2}
\end{gather*}
respectively, where
\begin{gather*}
f^2_A=g_B g_C,
\qquad
f^2_B=g_C g_A,
\qquad
f^2_C=g_A g_B,
\qquad
f_A f_B f_C =g_A g_B g_C.
\end{gather*}
The resulting equations assert that each of
\begin{gather}
 A+ \frac{qBC-q^{-1}CB}{q^2-q^{-2}},
\qquad B+
\frac{qCA-q^{-1}AC}{q^2-q^{-2}},
\qquad
C+
\frac{qAB-q^{-1}BA}{q^2-q^{-2}}
\label{eq:comlistpre}
\end{gather}
is a scalar in the underlying f\/ield.
We have eliminated
$g_A$, $g_B$, $g_C$ and are now
down to the three scalar parameters
(\ref{eq:comlistpre}).
To eliminate these we reinterpret them
as central elements in the algebra generated by
$A$, $B$, $C$. The resulting algebra is denoted $\Delta$
and called the
 universal Askey--Wilson algebra~\cite{uaw}.
The formal def\/inition of~$\Delta$
is given in Def\/inition~\ref{def:uaw} below.

In \cite{uaw}
we investigated~$\Delta$
from a ring theoretic point of view.
Our results include the following.
We displayed a faithful action of
the modular group ${\rm  {PSL}}_2(\mathbb Z)$ on $\Delta$
as a group of automorphisms
\cite[Theorems~3.1,~3.13]{uaw}.
We found several linear bases for~$\Delta$
\cite[Theorems~4.1,~7.5]{uaw}.
We described the center~$Z(\Delta)$
under the assumption that $q$ is not a root
of unity. For such~$q$ we found that
$Z(\Delta)$ is generated by the three
central elements~(\ref{eq:comlistpre}) together with an
element $\Omega$ called the Casimir element
\cite[Corollary~8.3]{uaw}.

We now discuss the equitable presentation for
$U_q(\mathfrak{sl}_2)$ \cite{equit}.
 This presentation involves generators
$x$, $y^{\pm 1}$, $z$
and relations $y y^{-1}=y^{-1}y = 1$,
\begin{gather*}
\frac{qxy-q^{-1}yx}{q-q^{-1}} = 1,
\qquad
\frac{qyz-q^{-1}zy}{q-q^{-1}} = 1,
\qquad
\frac{qzx-q^{-1}xz}{q-q^{-1}} = 1.
\end{gather*}
In \cite{alnajjar}
H.~AlNajjar investigated
Leonard pairs using the equitable
presentation of
$U_q(\mathfrak{sl}_2)$.
His approach is summarized as follows.
Let $V$ denote a f\/inite-dimensional irreducible
$U_q(\mathfrak{sl}_2)$-module.
Let $A$ denote an arbitrary
linear combination of
 $1$, $x$, $y$, $xy$ and let $B$ denote an arbitrary
linear combination of $1$, $y$, $z$, $yz$.
Consider the coef\/f\/icients.
Alnajjar
found necessary and suf\/f\/icient conditions on the coef\/f\/icients
for $A$, $B$ to
act on $V$
as a Leonard pair
\cite[Theorem~6.2]{alnajjar}.
In \cite{alnajjar2} Alnajjar  described the class of
Leonard pairs that result from his construction.
He showed that this class corresponds to
a family of orthogonal polynomials consisting of
the $q$-Racah, $q$-Hahn, dual $q$-Hahn, $q$-Krawtchouk,
dual $q$-Krawtchouk,
af\/f\/ine $q$-Krawtchouk, and quantum
$q$-Krawtchouk
polynomials.
For the Leonard pairs $A,B$ in the above class
 consider the corresponding Askey--Wilson relations.
We use the $\mathbb Z_3$-symmetric version in view of
the
$\mathbb Z_3$-symmetric nature of the
equitable presentation.
In the style of
Wiegmann and Zabrodin let
$C$ denote an arbitrary linear combination
of $1$, $z$, $x$, $zx$ and consider when $A$, $B$, $C$ satisfy some
$\mathbb Z_3$-symmetric Askey--Wilson relations.
Extending the work of Alnajjar one f\/inds that
 the ``most general'' solution is described as follows.

\begin{Proposition}
\label{prop:motiv}
Let $\mathbb F$ denote a field, and
fix a nonzero $q \in \mathbb F$ such that $q^4\not=1$.
Consider the $\mathbb F$-algebra
$U_q(\mathfrak{sl}_2)$
with equitable generators
$x$, $y^{\pm 1}$, $z$.
Let $a$, $b$, $c$ denote nonzero scalars in $\mathbb F$  and define
\begin{gather*}
A =  x a + y a^{-1} + \frac{xy-yx}{q-q^{-1}}bc^{-1},
\\
B =  y b + z b^{-1} + \frac{yz-zy}{q-q^{-1}}ca^{-1},
\\
C =  z c + x c^{-1} + \frac{zx-xz}{q-q^{-1}}ab^{-1}.
\end{gather*}
Then
\begin{gather*}
A+ \frac{qBC-q^{-1}CB}{q^2-q^{-2}}
=
\frac{
\Lambda (a+a^{-1})+(b+b^{-1})(c+c^{-1})
}{q+q^{-1}},
\\
B+
\frac{qCA-q^{-1}AC}{q^2-q^{-2}}
=
\frac{
\Lambda (b+b^{-1})+(c+c^{-1})(a+a^{-1})
}
{q+q^{-1}},
\\
C+
\frac{qAB-q^{-1}BA}{q^2-q^{-2}}
=
\frac{
\Lambda (c+c^{-1})+(a+a^{-1})(b+b^{-1})
}{q+q^{-1}}.
\end{gather*}
Here $\Lambda$ denotes the normalized
Casimir element of
$U_q(\mathfrak{sl}_2)$
from Lemma~{\rm \ref{lem:sixforms}} below.
\end{Proposition}

Let $q$, $a$, $b$, $c$ be from
Proposition~\ref{prop:motiv}. By that proposition
and since $\Lambda$ is central in
$U_q(\mathfrak{sl}_2)$,
there exists an algebra homomorphism
$\Delta \to
U_q(\mathfrak{sl}_2)$
that acts on the $\Delta$-generators
$A$, $B$, $C$ in the following way. It sends
\begin{gather*}
A \mapsto  x a + y a^{-1} + \frac{xy-yx}{q-q^{-1}}bc^{-1},
\\
B  \mapsto   y b + z b^{-1} + \frac{yz-zy}{q-q^{-1}}ca^{-1},
\\
C  \mapsto   z c + x c^{-1} + \frac{zx-xz}{q-q^{-1}}ab^{-1}.
\end{gather*}
It turns out that this homomorphism is not injective. In order
to shrink
the kernel
we reinterpret $a$, $b$, $c$ as mutually commuting indeterminates,
and view the above construction as giving an algebra homomorphism
$\Delta \mapsto
U_q(\mathfrak{sl}_2) \otimes_\mathbb F
\mathbb F \lbrack a^{\pm 1}, b^{\pm 1}, c^{\pm 1}\rbrack$.
A main result of the present paper is that this homomorphism is injective.
In another main result we compute the image of~$\Omega $ under the injection.
We also use the injection to show that~$\Delta$
contains no zero divisors.

The paper is organized as follows.
In Section~\ref{section2} we recall some
basic facts and then
state our main results, which are
Theorems \ref{thm:main}--\ref{thm:main3} and
Corollary
\ref{cor:nozerodiv}.
In Section~\ref{section3} we establish some identities
involving the equitable generators of
$U_q(\mathfrak{sl}_2)$,
which will be used
repeatedly.
In Section~\ref{section4} we prove Theorem~\ref{thm:main}.
In Section~\ref{section5} we prove Theorem~\ref{thm:main2}.
In Sections~\ref{section6}--\ref{section8}
we establish some slightly technical facts about
$U_q(\mathfrak{sl}_2)$,
which will be used in Section~\ref{section9}
to prove
Theorem~\ref{thm:main3}.
In Section~\ref{section10} we discuss some issues
concerning the
${\rm  {PSL}}_2(\mathbb Z)$ action on $\Delta$ that we mentioned earlier.

Our proofs for
Theorems~\ref{thm:main}--\ref{thm:main3}
are essentially self contained
and do not assume
Proposition~\ref{prop:motiv}.
We remark that
Proposition~\ref{prop:motiv} follows from Theorem~\ref{thm:main}.

For the rest of this paper $a$, $b$, $c$ denote mutually
commuting indeterminates.

\section{Statement of results}\label{section2}

Our conventions for the paper are as follows.
An algebra is meant to be associative and have a~1.
A subalgebra has the same~1 as the parent algebra.
We f\/ix a f\/ield~$\mathbb F$. All unadorned tensor products
are meant to be over $\mathbb F$. We f\/ix $q \in \mathbb F$ such that
$q^4\not=1$. Recall the natural numbers
$\mathbb N = \lbrace 0,1,2,\ldots \rbrace$ and
integers $\mathbb Z=\lbrace 0, \pm 1, \pm 2,\ldots \rbrace$.

\begin{Definition}[\protect{\cite[Def\/inition~1.2]{uaw}}]
\label{def:uaw}
Def\/ine an $\mathbb F$-algebra
$\Delta=\Delta_q$
by generators and relations in the following way.
The generators are $A$, $B$, $C$.
The relations assert that each of
\begin{gather}
 A+ \frac{qBC-q^{-1}CB}{q^2-q^{-2}},
\qquad B+
\frac{qCA-q^{-1}AC}{q^2-q^{-2}},
\qquad
C+
\frac{qAB-q^{-1}BA}{q^2-q^{-2}}
\label{eq:comlist}
\end{gather}
is central in $\Delta$.
The algebra
$\Delta$ is called the {\it universal Askey--Wilson algebra}.
\end{Definition}

\begin{Definition}[\protect{\cite[Def\/inition~1.3]{uaw}}]
\label{def:abc}
For the
three central elements in
(\ref{eq:comlist}), multiply each by $q+q^{-1}$ to get
$\alpha$, $\beta$, $\gamma$. Thus
\begin{gather}
A+ \frac{qBC-q^{-1}CB}{q^2-q^{-2}}
 =
\frac{\alpha}{q+q^{-1}},
\label{eq:u1}
\\
B+
\frac{qCA-q^{-1}AC}{q^2-q^{-2}}
 =
\frac{\beta}{q+q^{-1}},
\label{eq:u2}
\\
C+
\frac{qAB-q^{-1}BA}{q^2-q^{-2}}
 =
\frac{\gamma}{q+q^{-1}}.
\label{eq:u3}
\end{gather}
Note that each of $\alpha$, $\beta$, $\gamma$ is central in~$\Delta$.
\end{Definition}

We mention a few facts about~$\Delta$.
Recall that the modular group
${\rm  {PSL}}_2(\mathbb Z)$ has a presentation
by generators $\rho$, $\sigma$ and relations
$\rho^3=1$, $\sigma^2=1$. See for example
\cite{alpern}.
By \cite[Theorem~3.1]{uaw},
the group
${\rm  {PSL}}_2(\mathbb Z)$ acts on
$\Delta$ as a group of automorphisms
such that
$\rho$
 sends
$(A,B,C)\mapsto (B,C,A)$
and  $\sigma $
sends
$(A,B,\gamma)\mapsto (B,A,\gamma)$.
By
\cite[Theorem~3.13]{uaw}
this action is faithful.

By \cite[Theorem~4.1]{uaw} the following
is a basis for the $\mathbb F$-vector space $\Delta$:
\begin{gather*}
 A^iB^jC^k \alpha^r\beta^s\gamma^t,
\qquad
i,j,k,r,s,t \in \mathbb N.
\end{gather*}
There is a related basis
\cite[Theorem~7.5]{uaw} that we will use
 in
Section~\ref{section9} below. This related basis
involves
 a
central element $\Omega$ known as the Casimir
element~\cite[Lemma~6.1]{uaw}.
This element is def\/ined
as follows.

\begin{Definition}[\protect{\cite[Lemma~6.1]{uaw}}]
\label{def:casdelta}
Def\/ine $\Omega \in \Delta$ by
\begin{gather*}
\Omega = qABC+q^2A^2+q^{-2}B^2+q^2 C^2 - q A\alpha  - q^{-1} B \beta
- q C \gamma.
\end{gather*}
We call $\Omega$ the {\it Casimir element} of~$\Delta$.
\end{Definition}

\begin{Lemma}[\protect{\cite[Theorem~6.2, Corollary~8.3]{uaw}}]
The Casimir element $\Omega$ is contained in the
center~$Z(\Delta)$.
Moreover
$\lbrace \Omega^i\alpha^r\beta^s\gamma^t \,|\, i,r,s,t \in \mathbb N\rbrace$
is a basis for
the $\mathbb F$-vector space
$Z(\Delta)$,
provided that~$q$ is not a root of unity.
\end{Lemma}

\begin{Lemma}[\protect{\cite[Theorem~6.4]{uaw}}]
\label{lem:psl}
The Casimir element $\Omega$  is fixed by
everything in
${\rm  {PSL}}_2(\mathbb Z)$.
\end{Lemma}

We will be discussing how
$\Delta$ is related to the
quantum universal enveloping
algebra $U_q(\mathfrak{sl}_2)$.
For
this algebra
there are two presentations of
 interest to us;
 the Chevalley presentation
\cite[Section~1.1]{jantzen}
and the equitable presentation~\cite{equit}. We now recall the
  Chevalley presentation.

\begin{Definition}[\protect{\cite[Section~1.1]{jantzen}}]
\label{def:chev}
The
$\mathbb F$-algebra $U=U_q(\mathfrak{sl}_2)$
is def\/ined by generators
$e$, $f$, $k^{\pm 1 }$  and relations
\begin{gather*}
  k k^{-1} = k^{-1} k = 1,
\qquad
 ke = q^2 e k,
\qquad
 kf = q^{-2} f k,
\qquad ef-fe = \frac{k-k^{-1}}{q-q^{-1}}.
\end{gather*}
We call  $e$, $f$, $k^{\pm 1}$ the {\it Chevalley generators} for~$U$.
\end{Definition}

 We now brief\/ly discuss some
f\/inite-dimensional $U$-modules. Strictly speaking
we will not use this information; it is included
in order to clarify the nature of the Casimir element
for $U$ described below.

Recall the notation
\begin{gather*}
  \lbrack  n \rbrack_q=\frac{q^n-q^{-n}}{q-q^{-1}},
  \qquad n \in \mathbb N.
\end{gather*}

\begin{Lemma}[\protect{\cite[Section~2]{jantzen}}]
For all integers $n\geq 0$ and $\varepsilon \in \lbrace 1,-1\rbrace$
there exists a
$U$-module
$L(n,\varepsilon)$ with the following properties.
$L(n,\varepsilon)$ has a basis
$\lbrace v_i\rbrace_{i=0}^n$
such that
\begin{gather*}
k v_i  =  \varepsilon q^{n-2i} v_i, \quad  0 \leq i \leq n,
\\
f v_i = \lbrack i+1 \rbrack_q v_{i+1}, \quad 0 \leq i \leq n-1,
\qquad fv_n =0,
\\
e v_i = \varepsilon
\lbrack n-i+1 \rbrack_q v_{i-1}, \quad 1 \leq i \leq n,
\qquad ev_0 =0.
\end{gather*}
The
$U$-module
$L(n,\varepsilon)$ is irreducible provided that~$q$ is not a root of unity.
\end{Lemma}

In Def\/inition~\ref{def:casdelta} we gave the Casimir element
for $\Delta$. We now recall the Casimir element
for~$U$.

\begin{Definition}[\protect{\cite[Section~2.7]{jantzen}}]
\label{def:casu}
Def\/ine
 $\Phi \in
U$ as follows:
\begin{gather*}
\Phi = ef +\frac{q^{-1}k + qk^{-1}}{(q-q^{-1})^2}.
\end{gather*}
We call $\Phi$ the {\it Casimir element} of
$U$.
\end{Definition}

\begin{Lemma}[\protect{\cite[Lemma~2.7, Proposition~2.18]{jantzen}}]
The element
$\Phi$ is contained in the center~$Z(U)$.
Moreover $\lbrace \Phi^i\rbrace_{i \in \mathbb N}$ is a basis
for the $\mathbb F$-vector space~$Z(U)$,
 provided that $q$ is not a root of unity.
\end{Lemma}

\begin{Lemma}
\label{lem:normmeaning}
{\rm \cite[Lemma~2.7]{jantzen}.}
On the
$U$-module
$L(n,\varepsilon)$,
\begin{gather*}
\Phi  =  \varepsilon \frac{q^{n+1}+q^{-n-1}}{(q-q^{-1})^2} I.
\end{gather*}
Here $I$ denotes the identity map.
\end{Lemma}

For notational convenience we now adjust the
normalization for $\Phi$.

\begin{Definition}
\label{def:norm}
Def\/ine
\begin{gather}
\Lambda  = \big(q-q^{-1}\big)^2 \Phi
= \big(q-q^{-1}\big)^2ef +q^{-1}k + qk^{-1}.
\label{eq:lambdaform}
\end{gather}
Note that on
$L(n,\varepsilon)$,
\begin{gather*}
\Lambda  =  \varepsilon\big(q^{n+1}+q^{-n-1}\big)I.
\end{gather*}
We call
$\Lambda$ the {\it normalized Casimir element} for
$U$.
\end{Definition}

We now recall the equitable presentation for
$U$ \cite{equit}.

\begin{Proposition}[\protect{\cite[Theorem~2.1]{equit}}]
\label{prop:equit}
The algebra
$U$ is isomorphic to the
$\mathbb F$-algebra defined by generators
$x$, $y^{\pm 1}$, $z$ and relations
\begin{gather}
 y y^{-1} = y^{-1} y  = 1,
\label{eq:eq0}
\\
 \frac{qxy-q^{-1}yx}{q-q^{-1}} = 1,
\label{eq:eq1}
\\
 \frac{qyz-q^{-1}zy}{q-q^{-1}} = 1,
\label{eq:eq2}
\\
 \frac{qzx-q^{-1}xz}{q-q^{-1}} = 1.
\label{eq:eq3}
\end{gather}
An isomorphism with the presentation in
Definition~{\rm \ref{def:chev}} is given by
\begin{gather*}
 y^{\pm 1} \mapsto k^{\pm 1},
\qquad
 z \mapsto k^{-1} + f\big(q-q^{-1}\big),
\qquad
 x \mapsto k^{-1} - ek^{-1}q^{-1}\big(q-q^{-1}\big).
\end{gather*}
The inverse of this isomorphism is given by
\begin{gather*}
 k^{\pm 1} \mapsto y^{\pm 1},
\qquad
 f\mapsto (z-y^{-1})\big(q-q^{-1}\big)^{-1},
\qquad
 e\mapsto (1-xy)q\big(q-q^{-1}\big)^{-1}.
\end{gather*}
\end{Proposition}

\begin{Definition}[\protect{\cite[Def\/inition~2.2]{equit}}]
By the {\it equitable presentation} of~$U$ we mean the
presentation given in Proposition
\ref{prop:equit}. We call $x$, $y^{\pm 1}$, $z$ the
{\it equitable generators} for~$U$.
\end{Definition}

\begin{Note}
\label{note:ident}
In what follows we identify the copy of
$U$ given in
Def\/inition~\ref{def:chev}
with
the copy given in
Proposition~\ref{prop:equit},
via the isomorphism given in
Proposition~\ref{prop:equit}.
\end{Note}

In the equitable presentation
of $U$
the normalized Casimir element $\Lambda$ looks as follows.

\begin{Lemma}
\label{lem:sixforms}
The normalized Casimir element $\Lambda$ is equal to each of
the following:
\begin{alignat}{3}
& qx+q^{-1}y+qz-qxyz,
\qquad &&
q^{-1}x+qy+q^{-1}z-q^{-1}zyx, &
\label{eq:v1}
\\
& qy+q^{-1}z+qx-qyzx,
\qquad &&
q^{-1}y+qz+q^{-1}x-q^{-1}xzy,&
\label{eq:v2}
\\
& qz+q^{-1}x+qy-qzxy,
\qquad &&
q^{-1}z+qx+q^{-1}y-q^{-1}yxz.&
\label{eq:v3}
\end{alignat}
\end{Lemma}

\begin{proof}
For the data
(\ref{eq:v1})--(\ref{eq:v3})
let
$\Lambda^-_y$, $\Lambda^-_z$, $\Lambda^-_x$
denote the
expressions in the f\/irst column
and let
$\Lambda^+_y$, $\Lambda^+_z$, $\Lambda^+_x$ denote
the
expressions in the second column.
Consider the expression for
$\Lambda$ given in~(\ref{eq:lambdaform}).
Writing this expression in terms of $x$, $y$, $z$ using
the isomorphism in
Proposition~\ref{prop:equit} and Note~\ref{note:ident},
we obtain
$\Lambda=\Lambda^-_y$.
The element $\Lambda^-_y-\Lambda^+_z$ is equal to
$(q-q^{-1})x$ times
\begin{gather}
1 - \frac{qyz-q^{-1}zy}{q-q^{-1}}.
\label{eq:zer1}
\end{gather}
The expression
(\ref{eq:zer1}) is zero by
(\ref{eq:eq2}) so
$\Lambda^-_y=\Lambda^+_z$.
Similarly one f\/inds
$\Lambda^-_z=\Lambda^+_x$ and
$\Lambda^-_x=\Lambda^+_y$.
The element $\Lambda^-_y-\Lambda^+_x$ is equal to
\begin{gather}
1 - \frac{qxy-q^{-1}yx}{q-q^{-1}}
\label{eq:zer2}
\end{gather}
times
$(q-q^{-1})z$.
The expression~(\ref{eq:zer2}) is zero by~(\ref{eq:eq1}) so
$\Lambda^-_y=\Lambda^+_x$.
Similarly one f\/inds
$\Lambda^-_z=\Lambda^+_y$ and
$\Lambda^-_x=\Lambda^+_z$.
By these comments $\Lambda$ is equal to each of
$\Lambda^{\pm }_x$, $\Lambda^{\pm }_y$, $\Lambda^{\pm }_z$.
\end{proof}

Recall that $a$, $b$, $c$ are mutually commuting indeterminates.
Let
$\mathbb F \lbrack a^{\pm 1}, b^{\pm 1}, c^{\pm 1}\rbrack$
denote the $\mathbb F$-algebra
of Laurent polynomials in
$a$, $b$, $c$ that have all coef\/f\/icients in~$\mathbb F$.

We now state our main results.

\begin{Theorem}\label{thm:main}
There exists a unique  $\mathbb F$-algebra homomorphism
$\natural: \Delta \to
U \otimes
\mathbb F \lbrack a^{\pm 1}, b^{\pm 1}, c^{\pm 1}\rbrack
$ that sends
\begin{gather*}
A  \mapsto   x \otimes a + y \otimes a^{-1} +
\frac{x y - y x}{q-q^{-1}} \otimes b c^{-1},
\\
B \mapsto  y \otimes b + z \otimes b^{-1} +
\frac{y z - z y}{q-q^{-1}} \otimes c a^{-1},
\\
C \mapsto  z \otimes c + x \otimes c^{-1} +
\frac{z x - x z}{q-q^{-1}} \otimes a b^{-1},
\end{gather*}
where $x$, $y$, $z$ denote the equitable generators for~$U$.
The homomorphism $\natural$ sends
\begin{gather}
\alpha  \mapsto  \Lambda \otimes \big(a+a^{-1}\big) + 1\otimes \big(b+b^{-1}\big)\big(c+c^{-1}\big),
\label{eq:alsend}
\\
\beta  \mapsto  \Lambda \otimes \big(b+b^{-1}\big) + 1 \otimes \big(c+c^{-1}\big)\big(a+a^{-1}\big),
\label{eq:besend}
\\
\gamma  \mapsto  \Lambda \otimes \big(c+c^{-1}\big) + 1 \otimes \big(a+a^{-1}\big)\big(b+b^{-1}\big),
\label{eq:gasend}
\end{gather}
where $\Lambda $ denotes the normalized Casimir element of~$U$.
\end{Theorem}

\begin{Theorem}\label{thm:main2}
Under the homomorphism
$\natural$
from Theorem~{\rm \ref{thm:main}}, the image of
$\Omega$ is
\begin{gather}
1 \otimes \big(q+q^{-1}\big)^2
-
1\otimes \big(a+a^{-1}\big)^2
-
1\otimes \big(b+b^{-1}\big)^2
-
1 \otimes \big(c+c^{-1}\big)^2
\nonumber\\
\phantom{1 \otimes \big(q+q^{-1}\big)^2}{}
-\Lambda \otimes \big(a+a^{-1}\big)\big(b+b^{-1}\big)\big(c+c^{-1}\big)
- \Lambda ^2 \otimes 1.\label{eq:omimage}
\end{gather}
Here $\Lambda$ denotes the normalized Casimir element
of~$U$.
\end{Theorem}

\begin{Theorem}
\label{thm:main3}
The
homomorphism~$\natural$
from Theorem~{\rm \ref{thm:main}} is injective.
\end{Theorem}

We mention a corollary to
Theorem~\ref{thm:main3}.
For an $\mathbb F$-algebra $\mathcal A$,
an element $u \in {\mathcal A}$ is called a~{\it zero divisor}
whenever $u\not=0$ and
there exists $0 \not=v \in \mathcal A$
such that~$uv=0$.
By~\cite[Proposition~1.8]{jantzen}
the
algebra $U$
contains no zero divisors.
For an $\mathbb F$-algebra $\mathcal A$
and indeterminate $\lambda$ consider the $\mathbb F$-algebra
${\mathcal A}\otimes \mathbb F \lbrack \lambda,\lambda^{-1}\rbrack$.
One checks that
$\mathcal A$ contains no zero divisors if and only if
${\mathcal A}\otimes \mathbb F \lbrack \lambda,\lambda^{-1}\rbrack$
contains no zero divisors.
Applying this comment three times we see that
the algebra
$
U \otimes
\mathbb F \lbrack a^{\pm 1}, b^{\pm 1}, c^{\pm 1}\rbrack
$ contains no zero divisors.
By this and
Theorem~\ref{thm:main3} we obtain the following result.

\begin{Corollary}
\label{cor:nozerodiv}
The $\mathbb F$-algebra $\Delta$ contains no zero divisors.
\end{Corollary}

\section[The elements $\nu_x$, $\nu_y$, $\nu_z$]{The elements $\boldsymbol{\nu_x}$, $\boldsymbol{\nu_y}$, $\boldsymbol{\nu_z}$}\label{section3}

In this section we record a number of
identities involving
the equitable generators for
$U$.
These identities will be used in our
proof of Theorems
\ref{thm:main}--\ref{thm:main3}.

The relations
(\ref{eq:eq1})--(\ref{eq:eq3})
can be reformulated as follows:
\begin{gather*}
q(1-yz) = q^{-1}(1-zy),
\qquad
q(1-zx) = q^{-1}(1-xz),
\qquad
q(1-xy) = q^{-1}(1-yx).
\end{gather*}

\begin{Definition}
\label{def:nx}
Let $\nu_x$, $\nu_y$, $\nu_z$ denote the following elements
in
$U$:
\begin{gather}
 \nu_x = q(1-yz)=q^{-1}(1-zy),
\label{eq:defnx}
\\
 \nu_y = q(1-zx)=q^{-1}(1-xz),
\label{eq:defny}
\\
 \nu_z = q(1-xy)=q^{-1}(1-yx).
\label{eq:defnz}
\end{gather}
\end{Definition}

\begin{Note}
\label{note:efview}
We have
\begin{gather*}
 e =\frac{\nu_z}{q-q^{-1}},
\qquad
f =-\frac{q^{-1} y^{-1}\nu_x}{q-q^{-1}},
\qquad
 \nu_z = \big(q-q^{-1}\big)e, \qquad
\nu_x = -q\big(q-q^{-1}\big)kf.
\end{gather*}
\end{Note}

\begin{Lemma}
\label{lem:prod}
The following relations hold in~$U$:
\begin{alignat}{3}
& xy =1-q^{-1}\nu_z, \qquad && yx = 1-q\nu_z,&
\label{eq:xyx}
\\
& yz =1-q^{-1}\nu_x, \qquad &&  zy = 1-q\nu_x,&
\label{eq:yzy}
\\
& zx =1-q^{-1}\nu_y, \qquad && xz = 1-q\nu_y.&
\label{eq:zxz}
\end{alignat}
\end{Lemma}

\begin{proof}
These equations are
reformulations of
(\ref{eq:defnx})--(\ref{eq:defnz}).
\end{proof}

\begin{Lemma}
The following relations hold in
$U$:
\begin{alignat}{3}
& \frac{xy-yx}{q-q^{-1}}=
\nu_z, \qquad &&
\frac{qyx-q^{-1}xy}{q-q^{-1}}= 1-\big(q+q^{-1}\big)\nu_z,&
\label{eq:double1}
\\
& \frac{yz-zy}{q-q^{-1}}=
\nu_x,\qquad&&
\frac{qzy-q^{-1}yz}{q-q^{-1}}= 1-\big(q+q^{-1}\big)\nu_x,&
\label{eq:double2}
\\
& \frac{zx-xz}{q-q^{-1}}=
\nu_y,\qquad &&
\frac{qxz-q^{-1}zx}{q-q^{-1}}= 1-\big(q+q^{-1}\big)\nu_y. &
\label{eq:double3}
\end{alignat}
\end{Lemma}

\begin{proof}
For each equation
evaluate the left-hand side using
Lemma~\ref{lem:prod}.
\end{proof}

\begin{Lemma}
\label{lem:qcom}
The following relations hold in
$U$:
\begin{alignat}{3}
& x \nu_y = q^2 \nu_y x, \qquad &&
x \nu_z = q^{-2} \nu_z x, &
\label{eq:com1}
\\
& y \nu_z = q^2 \nu_z y, \qquad &&
y \nu_x = q^{-2} \nu_x y, &
\label{eq:com2}
\\
& z \nu_x = q^2 \nu_x z, \qquad &&
z \nu_y = q^{-2} \nu_y z. &
\label{eq:com3}
\end{alignat}
\end{Lemma}

\begin{proof}
Using $\nu_y = q(1-zx)$ we f\/ind
$q^{-1}x \nu_y=x-xzx$.
Using $\nu_y= q^{-1}(1-xz)$ we f\/ind
$q\nu_y x = x - xzx$. By these comments
$x \nu_y = q^2 \nu_y x$. The remaining
relations are similarly obtained.
\end{proof}

\begin{Lemma}
\label{lem:xvx}
The following relations hold in
$U$:
\begin{alignat}{3}
& \nu_x x = \Lambda - q y - q^{-1}z,
\qquad &&
x\nu_x = \Lambda - q^{-1}y - qz, &
\label{eq:nuxX}
\\
& \nu_y y = \Lambda - q z - q^{-1}x,
\qquad &&
y\nu_y = \Lambda - q^{-1}z - qx,&
\label{eq:nuyY}
\\
& \nu_z z = \Lambda - q x - q^{-1}y,
\qquad &&
z\nu_z = \Lambda - q^{-1}x - qy.&
\label{eq:nuzZ}
\end{alignat}
\end{Lemma}

\begin{proof}
To
verify the equation on
the left in
(\ref{eq:nuxX}), eliminate
$\nu_x$ using $\nu_x = q(1-yz)$, and eliminate
$\Lambda$
using the fact that $\Lambda$ is equal to
the expression on the left in~(\ref{eq:v2}).
 The remaining equations are
similarly verif\/ied.
\end{proof}

\begin{Lemma}
\label{lem:xnxcom}
The following relations hold in
$U$:
\begin{gather*}
\frac{x\nu_x-\nu_x x}{q-q^{-1}}  =  y-z,
\qquad
\frac{y\nu_y-\nu_y y}{q-q^{-1}}  =  z-x,
\qquad
\frac{z\nu_z-\nu_z z}{q-q^{-1}}  =  x-y.
\end{gather*}
\end{Lemma}

\begin{proof}
For each equation evaluate the left-hand side
using
Lemma~\ref{lem:xvx}.
\end{proof}

\begin{Lemma}
The normalized Casimir element $\Lambda$ is equal
to each of the following:
\begin{alignat}{3}
&\frac{qx \nu_x - q^{-1}\nu_x x}{q-q^{-1}}+\big(q+q^{-1}\big)z,
\qquad &&
\frac{q \nu_x x - q^{-1}x\nu_x}{q-q^{-1}}+\big(q+q^{-1}\big)y,&
\label{eq:dcas1}
\\
&\frac{qy \nu_y - q^{-1}\nu_y y}{q-q^{-1}}+\big(q+q^{-1}\big)x,
\qquad &&
\frac{q \nu_y y - q^{-1}y\nu_y}{q-q^{-1}}+\big(q+q^{-1}\big)z, &
\label{eq:dcas2}
\\
& \frac{qz \nu_z - q^{-1}\nu_z z}{q-q^{-1}}+\big(q+q^{-1}\big)y,
\qquad &&
\frac{q \nu_z z - q^{-1}z\nu_z}{q-q^{-1}}+\big(q+q^{-1}\big)x.&
\label{eq:dcas3}
\end{alignat}
\end{Lemma}

\begin{proof}
Evaluate each of the displayed expressions
using
Lemma~\ref{lem:xvx}.
\end{proof}

\begin{Lemma}
\label{lem:nxny}
The following relations hold in
$U$:
\begin{alignat}{3}
& \nu_x \nu_y = 1- q^{-1}\Lambda z + q^{-2} z^2,
\qquad &&
\nu_y \nu_x = 1- q\Lambda z + q^2 z^2,&
\label{eq:nxny1}
\\
& \nu_y \nu_z = 1- q^{-1}\Lambda x + q^{-2} x^2,
\qquad &&
\nu_z \nu_y = 1- q\Lambda x + q^2 x^2,&
\label{eq:nxny2}
\\
& \nu_z \nu_x = 1- q^{-1}\Lambda y + q^{-2} y^2,
\qquad &&
\nu_x \nu_z = 1- q\Lambda y + q^2 y^2.&
\label{eq:nxny3}
\end{alignat}
\end{Lemma}

\begin{proof}
To get the equation on the left in~(\ref{eq:nxny1}), observe
\begin{gather*}
\nu_x \nu_y  =  q^{-1}\nu_x (1-xz) =q^{-1}\nu_x -q^{-1} \nu_x x z
= 1-yz -q^{-1}\big(\Lambda - qy-q^{-1}z\big)z
\\
\phantom{\nu_x \nu_y}{}  =  1-q^{-1}\Lambda z + q^{-2}z^2.
\end{gather*}
The remaining equations are similarly verif\/ied.
\end{proof}

\begin{Lemma}
\label{lem:comnxny}
The following relations hold in
$U$:
\begin{gather}
\frac{q\nu_x\nu_y - q^{-1} \nu_y \nu_x}{q-q^{-1}}  =  1-z^2,
\label{eq:c1}
\\
\frac{q\nu_y\nu_z - q^{-1} \nu_z \nu_y}{q-q^{-1}}  =  1-x^2,
\label{eq:c2}
\\
\frac{q\nu_z\nu_x - q^{-1} \nu_x \nu_z}{q-q^{-1}}  =  1-y^2.
\label{eq:c3}
\end{gather}
\end{Lemma}

\begin{proof}
 For each equation evaluate the left-hand side using
Lemma~\ref{lem:nxny}.
\end{proof}

\section{The proof of Theorem  \ref{thm:main}}\label{section4}

In this section we prove Theorem~\ref{thm:main}.

For notational convenience we def\/ine some elements in
$U \otimes
\mathbb F \lbrack a^{\pm 1}, b^{\pm 1}, c^{\pm 1}\rbrack$:
\begin{gather}
A^\natural  =  x \otimes a + y \otimes a^{-1} + \nu_z \otimes bc^{-1},
\label{eq:aform}
\\
B^\natural  =  y \otimes b + z \otimes b^{-1} + \nu_x \otimes ca^{-1},
\label{eq:bform}
\\
C^\natural  =  z \otimes c + x \otimes c^{-1} + \nu_y \otimes ab^{-1}
\label{eq:cform}
\end{gather}
and
\begin{gather}
\alpha^\natural
 =
\Lambda \otimes \big(a+a^{-1}\big) + 1 \otimes \big(b+b^{-1}\big)\big(c+c^{-1}\big),
\label{eq:amaincom1}
\\
\beta^\natural  =
\Lambda \otimes \big(b+b^{-1}\big) + 1 \otimes \big(c+c^{-1}\big)(a+a^{-1}\big),
\label{eq:bmaincom2}
\\
\gamma^\natural
 =
\Lambda \otimes \big(c+c^{-1}\big) + 1 \otimes \big(a+a^{-1}\big)\big(b+b^{-1}\big).
\label{eq:cmaincom3}
\end{gather}
Note that each of $\alpha^\natural$, $\beta^\natural$, $\gamma^\natural$ is
central in
$U \otimes
\mathbb F \lbrack a^{\pm 1}, b^{\pm 1}, c^{\pm 1}\rbrack$.

\begin{proof}[Proof of Theorem~\ref{thm:main}.]
We f\/irst establish the existence of the homomorphism
in the theorem statement.
To do this it suf\/f\/ices to show that
\begin{gather}
A^\natural + \frac{qB^\natural C^\natural-q^{-1}C^\natural B^\natural}{q^2-q^{-2}}
 =
\frac{\alpha^\natural}{q+q^{-1}},
\label{eq:maincom1}
\\
B^\natural + \frac{qC^\natural A^\natural-q^{-1}A^\natural C^\natural}{q^2-q^{-2}}
 =
\frac{\beta^\natural}{q+q^{-1}},
\label{eq:maincom2}
\\
C^\natural + \frac{qA^\natural B^\natural-q^{-1}B^\natural A^\natural}{q^2-q^{-2}}
=
\frac{\gamma^\natural}{q+q^{-1}}.
\label{eq:maincom3}
\end{gather}
We
verify~(\ref{eq:maincom3}).
Let $P$ denote the
left-hand side of~(\ref{eq:maincom3})
minus
the right-hand side of
(\ref{eq:maincom3}).
We show that $P=0$.
View $P$ as a Laurent
polynomial in $a$, $b$, $c$ that has coef\/f\/icients in~$U$.
We will show that in this polynomial each coef\/f\/icient is zero.
To this end, evaluate
 $P$ using~(\ref{eq:aform})--(\ref{eq:cform}), (\ref{eq:cmaincom3}) and then collect terms.
We list below the terms for which the coef\/f\/icient in $P$ is potentially
nonzero:
\begin{gather}
\label{eq:poterms}
 ab,\quad ab^{-1}, \quad a^{-1}b,\quad a^{-1}b^{-1},\quad
 c, \quad c^{-1}, \quad a^{-2}c,\quad b^2c^{-1}.
\end{gather}
For each of these terms the coef\/f\/icient in $P$ is listed
in the table below, along with a
 reason why that
coef\/f\/icient
is zero.

\begin{center}
\begin{tabular}{c| c | c}
term  &  coef\/f\/icient in $P$ & why the coef\/f\/icient is 0
\\
\hline
\hline
$ab$ & $\frac{qxy-q^{-1}yx}{q^2-q^{-2}}-\frac{1}{q+q^{-1}}$\tsep{4pt}
& equation (\ref{eq:eq1})
\\
$ab^{-1}$
&
$\frac{qxz-q^{-1}zx}{q^2-q^{-2}}+\nu_y -\frac{1}{q+q^{-1}}$\tsep{4pt}
&
equation (\ref{eq:double3})
\\
$a^{-1}b$ &
$
\frac{q(y^2+ \nu_z\nu_x)-q^{-1}(y^2+\nu_x \nu_z)}{q^2-q^{-2}}\tsep{4pt}
-\frac{1}{q+q^{-1}}
$
&
equation (\ref{eq:c3})
\\
$a^{-1}b^{-1}$ &
$\frac{qyz-q^{-1}zy}{q^2-q^{-2}}-\frac{1}{q+q^{-1}}$\tsep{4pt}\bsep{4pt}
&
equation (\ref{eq:eq2})
\\
\hline
$c$
&
$
\frac{q x\nu_x-q^{-1}\nu_x x}{q^2-q^{-2}}+ z - \frac{\Lambda}{q+q^{-1}}
$\tsep{4pt}
&
equation
(\ref{eq:dcas1})
\\
$c^{-1}$ & $\frac{q\nu_z z - q^{-1} z \nu_z}{q^2-q^{-2}}+x -
\frac{\Lambda}{q+q^{-1}} $\tsep{4pt}\bsep{4pt}
&
equation (\ref{eq:dcas3})
\\
\hline
$a^{-2}c$
&
$\frac{q y \nu_x - q^{-1} \nu_x y}{q^2-q^{-2}} $\tsep{4pt}
&
equation (\ref{eq:com2})
\\
$b^2c^{-1}$
&
$\frac{q \nu_z y - q^{-1} y \nu_z}{q^2-q^{-2}} $\tsep{4pt}
&
equation (\ref{eq:com2})
\end{tabular}
        \end{center}

We have shown that for each term in
(\ref{eq:poterms}) the coef\/f\/icient in $P$ is zero.
Therefore $P=0$ and
the equation
(\ref{eq:maincom3}) holds.
The equations~(\ref{eq:maincom1}),
(\ref{eq:maincom2})
are similarly verif\/ied.
We have shown that the homomorphism exists. The homomorphism
is unique since $A$, $B$, $C$ generate $\Delta$.
The homomorphism satisf\/ies
(\ref{eq:alsend})--(\ref{eq:gasend}) by
(\ref{eq:amaincom1})--(\ref{eq:cmaincom3})
and
(\ref{eq:maincom1})--(\ref{eq:maincom3}).
\end{proof}

\section{The proof of Theorem~\ref{thm:main2}}\label{section5}

In this section we prove Theorem~\ref{thm:main2}.
Recall the Casimir element~$\Omega$ from Def\/inition~\ref{def:casdelta}.

\begin{proof}[Proof of Theorem~\ref{thm:main2}.]
By
Theorem
\ref{thm:main}
the image of $\Omega$ is
\begin{gather}
\label{eq:omegacheck}
qA^\natural B^\natural C^\natural +
q^2 (A^{\natural})^2 +
q^{-2} (B^{\natural})^2 +
q^{2} (C^{\natural})^2 -
q A^\natural \alpha^\natural -
q^{-1} B^\natural \beta^\natural -
q C^\natural\gamma^\natural,
\end{gather}
where
$A^\natural$, $B^\natural,C^\natural$,
$\alpha^\natural$, $\beta^\natural$,
$\gamma^\natural$ are from
(\ref{eq:aform})--(\ref{eq:cmaincom3}).
We show that~(\ref{eq:omegacheck})
is equal to~(\ref{eq:omimage}).
Def\/ine~$Q$ to be~(\ref{eq:omegacheck})
minus~(\ref{eq:omimage}).
We show
that $Q=0$. To do this we proceed as in the proof of
Theorem~\ref{thm:main}.
View $Q$ as a Laurent polynomial in
$a$, $b$, $c$ that has all coef\/f\/icients in~$U$.
 We will show that for this polynomial each
 coef\/f\/icient is zero.
To this end, evaluate~$Q$
using
(\ref{eq:aform})--(\ref{eq:cmaincom3}) and then collect terms.
Below we list the terms for which the coef\/f\/icient in~$Q$
is potentially nonzero:
\begin{gather}
a^2, \quad  a^{-2}, \quad a^2b^{-2}, \quad abc^{-1}, \quad ab^{-1}c^{-1},
\label{eq:round1}
\\
b^2, \quad b^{-2}, \quad b^2c^{-2}, \quad bca^{-1}, \quad bc^{-1}a^{-1},
\label{eq:round2}
\\
c^2, \quad c^{-2}, \quad  c^2a^{-2}, \quad cab^{-1},\quad  ca^{-1}b^{-1},
\label{eq:round3}
\\
abc, \quad a^{-1}b^{-1}c^{-1}, \quad 1.
\label{eq:round4}
\end{gather}
We show that for each term in
(\ref{eq:round1})--(\ref{eq:round4}) the coef\/f\/icient in~$Q$ is zero.
The coef\/f\/icient of $a^2$ in~$Q$~is
\begin{gather}
q x y \nu_y + q^2 x^2 -q x \Lambda -q \nu_y +1.
\label{eq:a2coef}
\end{gather}
To see that
(\ref{eq:a2coef})
is zero,  eliminate $xy$ using
the equation on the left in
(\ref{eq:xyx}), and evaluate
the result
using the equation on the right in
(\ref{eq:nxny2}).
The coef\/f\/icient of $a^{-2}$ in $Q$ is
\begin{gather}
\label{eq:am2coef}
qy \nu_x x +q^2 y^2-q \Lambda y -q^{-1} \nu_x +1.
\end{gather}
To see that
(\ref{eq:am2coef}) is zero,
f\/irst eliminate~$y \nu_x$ using
the equation on the right in~(\ref{eq:com2}). Evaluate the
result using the equation on the right in~(\ref{eq:xyx}) followed by the equation on the
right in~(\ref{eq:nxny3}).
The coef\/f\/icient of
$a^2b^{-2}$ in $Q$ is
\begin{gather}
\label{eq:a2bm2coef}
qx z \nu_y +q^2 \nu_y^2-q \nu_y.
\end{gather}
The expression~(\ref{eq:a2bm2coef}) is zero
by the equation on the right in~(\ref{eq:zxz}).
The coef\/f\/icient of~$abc^{-1}$ in~$Q$~is
\begin{gather}
\label{eq:abcmcoeff}
q(xyx + \nu_z y \nu_y)+ q^2(x \nu_z+\nu_z x)
-q(x+\Lambda \nu_z)
-q^{-1} y - q x + \Lambda.
\end{gather}
To see that~(\ref{eq:abcmcoeff})
is zero,
eliminate~$yx$ using the equation on the right
in~(\ref{eq:xyx}),
and eliminate~$y \nu_y$ using
the equation on the right in~(\ref{eq:nuyY}).
Simplify the result
using
the equation on the left in~(\ref{eq:nuzZ}).
The coef\/f\/icient of~$ab^{-1}c^{-1}$ in $Q$ is
\begin{gather}
\label{eq:abmcmcoeff}
q(xzx + \nu_z z \nu_y)+ q^2(x \nu_y+\nu_y x)
- q x
-q^{-1}z
-q(x+\Lambda \nu_y)
+ \Lambda.
\end{gather}
To see that
(\ref{eq:abmcmcoeff}) is zero,
eliminate
$xz$ using the equation on the right in~(\ref{eq:zxz}), and
eliminate~$\nu_z z$ using the
equation on the left in~(\ref{eq:nuzZ}).
Simplify the result
using
the equation on the right in~(\ref{eq:nuyY}).
We have shown that for each term in~(\ref{eq:round1}) the coef\/f\/icient in~$Q$ is zero.
By Lemma~\ref{lem:psl} $\Omega$ is f\/ixed by
the automorphism of~$\Delta$ that sends
$(A,B,C)$ to
$(B,C,A)$.
Combining this with the~$\mathbb Z_3$-symmetric nature
of
(\ref{eq:aform})--(\ref{eq:cmaincom3}), we see that
for each term in~(\ref{eq:round2}),~(\ref{eq:round3}) the coef\/f\/icient in~$Q$ is also zero.
We now consider the terms in~(\ref{eq:round4}).
The coef\/f\/icient of~$abc$ in~$Q$ is
\begin{gather}
qxyz -qx -q^{-1}y-qz+\Lambda.
\label{eq:abccoeff}
\end{gather}
The expression (\ref{eq:abccoeff}) is zero
using the left side of~(\ref{eq:v1}).
The coef\/f\/icient of $a^{-1}b^{-1}c^{-1}$ in $Q$ is
\begin{gather}
\label{eq:abcmmmcoeff}
qyzx-qy-q^{-1}z -qx +\Lambda.
\end{gather}
The expression~(\ref{eq:abcmmmcoeff}) is zero
using the left side of~(\ref{eq:v2}).
The constant term in $Q$ is
\begin{gather*}
 q\big(\nu_z z^2 + x \nu_x x + y^2 \nu_y + \nu_z \nu_x \nu_y\big)
+ q^2(xy + yx)
+ q^{-2}(yz+zy)
  +  q^2(zx+xz)\\
\qquad{} -q(\nu_z+ \Lambda y+\Lambda x)
-q^{-1}(\nu_x+\Lambda z+ \Lambda y)
 -q(\nu_y+\Lambda x+\Lambda z)
-\big(q+q^{-1}\big)^2+\Lambda^2+6.
\end{gather*}
We show that this constant term is equal to zero.
Using Lemma~\ref{lem:prod}
and Lemma~\ref{lem:xvx}
we f\/ind
\begin{gather}
\nu_z z^2  =  q^{-2}\nu_x + q^2\nu_y+\Lambda z  - q -q^{-1},
\label{eq:1}
\\
x \nu_x x  =  \nu_y + \nu_z+\Lambda x  - q -q^{-1},
\label{eq:2}
\\
y^2 \nu_y  =  q^{-2}\nu_x + q^2\nu_z+\Lambda y  - q -q^{-1}.
\label{eq:3}
\end{gather}
Using the equation on the left in
(\ref{eq:nxny1}), followed by
(\ref{eq:1}) and the equation on the left in
(\ref{eq:nuzZ}),
we f\/ind
\begin{gather}
\nu_z \nu_x \nu_y =
q^{-4}\nu_x
+
\nu_y
+
\nu_z
+\Lambda x + q^{-2}\Lambda y+q^{-2}\Lambda z-q^{-1}\Lambda^2-q^{-1}-q^{-3}.
\label{eq:4}
\end{gather}
By Lemma
\ref{lem:prod},
\begin{gather}
xy+yx = 2-\big(q+q^{-1}\big)\nu_z,
\\
yz+zy = 2-\big(q+q^{-1}\big)\nu_x,
\\
zx+xz = 2-\big(q+q^{-1}\big)\nu_y.
\label{eq:5}
\end{gather}
Simplifying the constant term of $Q$ using
(\ref{eq:1})--(\ref{eq:5}) we f\/ind that this constant term is equal
to~zero.
We have shown that for each term in
(\ref{eq:round1})--(\ref{eq:round4}) the coef\/f\/icient in~$Q$ is zero.
Therefore \mbox{$Q=0$} and the result follows.
\end{proof}

\section[A $\mathbb Z$-grading of $U$]{A $\boldsymbol{\mathbb Z}$-grading of $\boldsymbol{U}$}\label{section6}

Our next general goal is to prove Theorem~\ref{thm:main3}.
To prepare for this proof we obtain some results
about~$U$.
In this section we discuss
a certain
$\mathbb Z$-grading of~$U$.
In the next section we will use this
$\mathbb Z$-grading of~$U$
to
 get
a
$\mathbb Z$-grading of
$U \otimes
\mathbb F \lbrack a^{\pm 1}, b^{\pm 1}, c^{\pm 1}\rbrack $.
The
$\mathbb Z$-grading of
$U\otimes
\mathbb F \lbrack a^{\pm 1}, b^{\pm 1}, c^{\pm 1}\rbrack $
will be used in our proof of
Theorem~\ref{thm:main3}.

Let $\mathcal A$ denote an
$\mathbb F$-algebra.
By a {\it $\mathbb Z$-grading of $\mathcal A$} we mean a sequence
 $\lbrace {\mathcal A}_n \rbrace_{n\in \mathbb Z}$ consisting of
 subspaces of~$\mathcal A$
 such
that
\begin{gather*}
{\mathcal A} = \sum_{n\in \mathbb Z} \mathcal A_n
\qquad
{\mbox {\rm {(direct sum)}}},
\end{gather*}
and $\mathcal A_m \mathcal A_n \subseteq \mathcal A_{m+n}$
for all $m,n \in \mathbb Z$. Let
$\lbrace {\mathcal A}_n \rbrace_{n\in \mathbb Z}$ denote
a $\mathbb Z$-grading of~$\mathcal A$.
For $n \in \mathbb Z$
we call~$\mathcal A_n$ the {\it $n$-homogeneous component} of~$\mathcal A$.
We refer to~$n$ as the {\it degree}
of~$\mathcal A_n$.
An element of~$\mathcal A$ is
said to be {\it homogeneous with degree $n$} whenever
it is contained in~$\mathcal A_n$. Pick $\xi \in \mathcal A$
and write
$\xi = \sum\limits_{n \in \mathbb Z} \xi_n$ with
$\xi_n \in \mathcal A_n$ for $n \in \mathbb Z$.
We call the elements $\lbrace \xi_n\rbrace_{n\in \mathbb Z}$
the {\it homogeneous components of~$\xi$}.

\begin{Lemma}[\protect{\cite[Theorem~1.5]{jantzen}}]
\label{lem:efkbasis}
The following is a basis for the $\mathbb F$-vector space
$U$:
\begin{gather}
e^r k^s f^t, \qquad   r,t\in \mathbb N, \qquad  s \in \mathbb Z.
\label{eq:kef}
\end{gather}
\end{Lemma}

For $n \in \mathbb Z$ let $U_n$ denote the subspace of
$U$  spanned
by those elements
$e^r k^s f^t$ from~(\ref{eq:kef}) that satisfy
$r-t=n$.
By \cite[Section~1.9]{jantzen} the sequence
$\lbrace U_n\rbrace_{n\in \mathbb Z}$ is a $\mathbb Z$-grading
of
$U$.
With respect to this $\mathbb Z$-grading the elements
$e$, $k$, $f$ are homogeneous with degrees
$1$, $0$, $-1$ respectively. Moreover the normalized
Casimir element $\Lambda$ from~\eqref{eq:lambdaform}
is homogeneous with degree~0.

By construction,
for $n \in \mathbb Z$ the $n$-homogeneous component
$U_n$ has a basis consisting of the elements
$e^r k^s f^t$ from~(\ref{eq:kef}) that satisfy
$r-t=n$. There is another basis for
$U_n$ that is better suited to our purpose;
this basis involves $\Lambda $ and
 will be displayed shortly.

\begin{Lemma}
\label{lem:ident}
For an integer $t\geq 0$,
\[
e^t f^t = \prod_{i=1}^t \frac{\Lambda -q^{1-2i}k-q^{2i-1}k^{-1}}{(q-q^{-1})^2}.
\]
\end{Lemma}

\begin{proof}
Assume $t\geq 1$; otherwise the result is trivial.
Using
(\ref{eq:lambdaform}) and $ke=q^2ek$,
\begin{gather*}
e^t f^t  =  e^{t-1} ef f^{t-1}
= e^{t-1} \frac{\Lambda -q^{-1}k-qk^{-1}}{(q-q^{-1})^2}f^{t-1}
 = \frac{\Lambda -q^{1-2t}k-q^{2t-1}k^{-1}}{(q-q^{-1})^2}e^{t-1}f^{t-1}.
\end{gather*}
The result follows by induction on~$t$.
\end{proof}

\begin{Lemma}
\label{lem:basisalt}
For all integers $n\geq 0$ the following $(i)$, $(ii)$ hold.
\begin{itemize}\itemsep=0pt
\item[$(i)$]
The $\mathbb F$-vector space $U_n$ has a basis
\begin{gather*}
e^n k^s \Lambda^t,
 \qquad s\in \mathbb Z, \qquad   t\in \mathbb N.
\end{gather*}
\item[$(ii)$]
The $\mathbb F$-vector space $U_{-n}$ has a basis
\begin{gather*}
k^s \Lambda^t  f^n,
 \qquad s\in \mathbb Z, \qquad   t\in \mathbb N.
\end{gather*}
\end{itemize}
\end{Lemma}

\begin{proof}
$(i)$ The elements
$\lbrace k^s\rbrace_{s \in \mathbb Z}$ are linearly
independent
by
Lemma~\ref{lem:efkbasis}, so they form
a basis for a subalgebra of~$U$ which we denote by~$K$.
By Lemma~\ref{lem:efkbasis}
the sum $U_n=\sum\limits_{\ell=0}^\infty e^{n+\ell} K f^\ell$ is direct.
We have
$Ke=eK$ since
$ke=q^2ek$, and
similarly $Kf=fK$.
Pick an integer $t\geq 0$.
By Lemma~\ref{lem:ident}  and induction on $t$
we f\/ind
$ \Lambda^t \in
\sum\limits_{\ell=0}^t e^\ell Kf^\ell$
and
$\Lambda^t -(q-q^{-1})^{2t} e^t f^t
\in
\sum\limits_{\ell=0}^{t-1} e^{\ell} K f^\ell $.
For the above $t$ and all $s\in \mathbb Z$
we have $e^n k^s \Lambda^t  \in
\sum\limits_{\ell=0}^t e^{n+\ell} K f^\ell $ and
\begin{gather*}
e^n k^s \Lambda^t - \big(q-q^{-1}\big)^{2t} q^{2st}e^{n+t}k^s f^t \in
\sum_{\ell=0}^{t-1} e^{n+\ell} K f^\ell.
\end{gather*}
The result follows
from these comments
and the fact that
$\lbrace
e^{n+t}k^s f^t \,|\, s\in \mathbb Z,\;t\in \mathbb N\rbrace$ is a basis for~$U_n$.

$(ii)$. Similar to the proof of $(i)$ above.
\end{proof}

 We now consider the $\mathbb Z$-grading
$\lbrace U_n\rbrace_{n\in \mathbb Z}$ from the point
of view of the equitable presentation.

\begin{Lemma}
\label{lem:gensetu}
The $\mathbb F$-algebra $U$ is generated by
$\nu_x$, $y^{\pm 1}$, $\nu_z$. Moreover
\begin{gather}
x = y^{-1}- q^{-1}\nu_z y^{-1},
\qquad
z = y^{-1}- q^{-1}y^{-1} \nu_x.
\label{eq:xzform}
\end{gather}
\end{Lemma}

\begin{proof}
The equation
on the left
in~(\ref{eq:xzform})
is a reformulation of
the equation on the left
in~(\ref{eq:xyx}).
The equation on the right in~(\ref{eq:xzform}) is similarly obtained.
The f\/irst assertion of the lemma follows
from~(\ref{eq:xzform}) and the fact that
$x$, $y^{\pm 1}$, $z$ generate
$U$.
\end{proof}

\begin{Lemma}
\label{lem:weightlist}
The generators
$\nu_x$, $y^{\pm 1}$, $\nu_z$ are homogeneous
with degree $-1$, $0$, $1$ respectively.
\end{Lemma}
\begin{proof}
Use Note
\ref{note:efview} and $y=k$, along with
 the comments below
Lemma~\ref{lem:efkbasis}.
\end{proof}

\begin{Lemma}
\label{lem:unbasis}
Pick an integer $n\geq 0$.
The $\mathbb F$-vector space $U_n$ has a basis
\begin{gather*}
\nu^n_z  y^i \Lambda^j,
  \qquad i\in \mathbb Z, \qquad j\in \mathbb N.
\end{gather*}
The $\mathbb F$-vector space $U_{-n}$ has a basis
\begin{gather*}
y^i \Lambda^j \nu^n_x,
 \qquad i\in \mathbb Z, \qquad j\in \mathbb N.
\end{gather*}
\end{Lemma}
\begin{proof}
This is a reformulation of
Lemma~\ref{lem:basisalt}, using
Note~\ref{note:efview} and $y=k$.
\end{proof}

 We comment on the homogeneous component
$U_0$.

\begin{Lemma}
\label{lem:u0facts}
The homogeneous component
$U_0$ is the $\mathbb F$-subalgebra
of~$U$ generated by
$y^{\pm 1}$,
$\Lambda $.
The algebra~$U_0$
is commutative.
The following is a basis for the
$\mathbb F$-vector space~$U_0$:
\begin{gather}
\label{eq:u0basis}
 y^i \Lambda^j,
  \qquad i\in \mathbb Z, \qquad    j\in \mathbb N.
\end{gather}
\end{Lemma}
\begin{proof}
To get the basis
(\ref{eq:u0basis}) set $n=0$ in
Lemma~\ref{lem:unbasis}. The remaining assertions
are clear.
\end{proof}

Let $\lambda_1$, $\lambda_2$ denote commuting
indeterminates.

\begin{Corollary}
\label{cor:laurant}
There exists an  $\mathbb F$-algebra isomorphism
$U_0 \to \mathbb F\lbrack \lambda_1, \lambda^{\pm 1}_2]$
that sends
$\Lambda \to \lambda_1$
and
$y \to \lambda_2$.
\end{Corollary}

\begin{proof}
Immediate from Lemma
\ref{lem:u0facts}.
\end{proof}

\begin{Definition}
\rm For $n \in \mathbb Z$ def\/ine an $\mathbb F$-linear map
$\pi_n:
U \to
U
$
such that $(\pi_n-1)U_n=0$ and
$\pi_n U_m=0$ if $m\not=n$  ($m\in \mathbb Z)$.
Thus $\pi_n$ is the projection
from $U $ onto $U_n$.
Note that for $u \in
U$ the element
$\pi_n(u)$
is the homogeneous component of $u$ with degree~$n$.
\end{Definition}

\begin{Lemma}
\label{lem:firsttable}
In the table below we list some
elements $u \in
U$.
For each element $u$ we display the homogeneous
component
$\pi_n(u)$ for $-1\leq n \leq 1$.
All other homogeneous components of $u$ are zero.
\begin{center}
\begin{tabular}{c| c  c c}
$u$ &  $\pi_{-1}(u)$ & $ \pi_0(u)$  & $ \pi_1(u)$
\\
\hline
\hline
$x$ & $0$ & $y^{-1}$ & $-q^{-1} \nu_z y^{-1}$\tsep{2pt}
\\
$y$ & $0$ & $y$ & $0$
\\
$z$ & $-q^{-1}y^{-1}\nu_x$ & $y^{-1}$ & $0$
\\
\hline
$\nu_x$ & $\nu_x$ & $0$ & $0$
\\
$\nu_y$ & $q^{-2}y^{-2}\nu_x$ & $y^{-1}\Lambda-(q+q^{-1})y^{-2}$
& $q^{-2}\nu_z y^{-2}$
\\
$\nu_z$ & $0$ & $0$ & $\nu_z$
\\
\hline
$\Lambda $ & $0$ & $\Lambda $ & $0$
\end{tabular}
        \end{center}
\end{Lemma}
\begin{proof}
The assertions about
$y$, $\nu_x$, $\nu_z$ come  from
Lemma~\ref{lem:weightlist}.
The assertions about $x$, $z$ follow from
(\ref{eq:xzform}).
We mentioned below
Lemma~\ref{lem:efkbasis}
that~$\Lambda $ is homogeneous with degree 0.
To verify the assertion about
$\nu_y$, in the equation $ \nu_y = y^{-1} y \nu_y$
evaluate~$y\nu_y $ using
the equation on the right in~(\ref{eq:nuyY}), and simplify the
result using
rows $x$, $z$ of the above table
along with the
equation on the left in~(\ref{eq:com2}).
\end{proof}

\section[A $\mathbb Z$-grading of $U\otimes \mathbb F \lbrack a^{\pm 1}, b^{\pm 1}, c^{\pm 1}\rbrack $]{A $\boldsymbol{\mathbb Z}$-grading of $\boldsymbol{U\otimes \mathbb F \lbrack a^{\pm 1}, b^{\pm 1}, c^{\pm 1}\rbrack}$}\label{section7}

In the previous section we discussed a
$\mathbb Z$-grading of~$U$.
In the present section we extend this
to a
$\mathbb Z$-grading of
$U \otimes
\mathbb F \lbrack a^{\pm 1}, b^{\pm 1}, c^{\pm 1}\rbrack $.
We obtain some results about
the
 $\mathbb Z$-grading
of
$U \otimes
\mathbb F \lbrack a^{\pm 1}, b^{\pm 1}, c^{\pm 1}\rbrack $
that will be used to prove
Theorem~\ref{thm:main3}.

By Lemma
\ref{lem:gensetu}
the elements $\nu_x$, $y^{\pm 1}$, $\nu_z$ form a generating
set for the $\mathbb F$-algebra~$U$.
Therefore the following is a generating set
for the $\mathbb F$-algebra
$U \otimes
\mathbb F \lbrack a^{\pm 1}, b^{\pm 1}, c^{\pm 1}\rbrack $:
\begin{gather}
\nu_x\otimes 1,\quad
y^{\pm 1}\otimes 1,\quad
\nu_z\otimes 1,\quad
1 \otimes a^{\pm 1},
\quad
1 \otimes b^{\pm 1},
\quad
1 \otimes c^{\pm 1}.
\label{eq:expg}
\end{gather}
Consider the $\mathbb Z$-grading
of
$U$
from below
Lemma~\ref{lem:efkbasis}.
This
$\mathbb Z$-grading
of
$U$
induces a
$\mathbb Z$-grading of
$U\otimes
\mathbb F \lbrack a^{\pm 1}, b^{\pm 1}, c^{\pm 1}\rbrack $
whose
homogeneous components are described as follows.
For $n \in \mathbb Z$ the $n$-homogeneous component
is $U_n \otimes
\mathbb F \lbrack a^{\pm 1}, b^{\pm 1}, c^{\pm 1}\rbrack $.
With respect to this $\mathbb Z$-grading of
$U \otimes
\mathbb F \lbrack a^{\pm 1}, b^{\pm 1}, c^{\pm 1}\rbrack $
the generators~(\ref{eq:expg}) are homogeneous with the following
degrees:

\begin{center}
\begin{tabular}{c||c c c |c c c}
$v$ &  $\nu_x \otimes 1$ & $y^{\pm 1} \otimes 1$  & $\nu_z\otimes 1$
&
$1 \otimes a^{\pm 1} $ &
$1 \otimes b^{\pm 1} $ &
$1 \otimes c^{\pm 1} $
\\
\hline
degree of $v$ & $-1$ & $0$ & $1$ & $0$ & $0$ & $0$
\end{tabular}
 \end{center}

\begin{Lemma}
\label{lem:gradingbasis}
Pick an integer $n\geq 0$. The $\mathbb F$-vector space
$U_n \otimes
\mathbb F \lbrack a^{\pm 1}, b^{\pm 1}, c^{\pm 1}\rbrack $ has
a basis
\begin{gather*}
\nu^n_z y^i \Lambda^j
\otimes a^rb^sc^t,
 \qquad i,r,s,t\in \mathbb Z, \qquad j\in \mathbb N.
\end{gather*}
The $\mathbb F$-vector space $U_{-n} \otimes {\mathbb F}\lbrack a^{\pm 1}, b^{\pm 1}, c^{\pm 1} \rbrack$ has a basis
\begin{gather*}
y^i \Lambda^j \nu^n_x
\otimes a^rb^sc^t,
 \qquad i,r,s,t\in \mathbb Z,  \qquad j\in \mathbb N.
\end{gather*}
\end{Lemma}
\begin{proof}
By Lemma
\ref{lem:unbasis}
and the construction.
\end{proof}

We comment on the homogeneous component
$U_0 \otimes
\mathbb F \lbrack a^{\pm 1}, b^{\pm 1}, c^{\pm 1}\rbrack $.

\begin{Lemma}
\label{lem:u0tensor}
The homogeneous component
$U_0 \otimes
\mathbb F \lbrack a^{\pm 1}, b^{\pm 1}, c^{\pm 1}\rbrack $
is the $\mathbb F$-subalgebra of
$
U
\otimes
\mathbb F \lbrack a^{\pm 1}, b^{\pm 1}, c^{\pm 1}\rbrack $
generated by
\begin{gather*}
y^{\pm 1}\otimes 1,\quad
\Lambda \otimes 1, \quad
1\otimes a^{\pm 1},
\quad
1\otimes b^{\pm 1},
\quad
1\otimes c^{\pm 1}.
\end{gather*}
The algebra
$U_0 \otimes
\mathbb F \lbrack a^{\pm 1}, b^{\pm 1}, c^{\pm 1}\rbrack $
is commutative.
The following is a basis for the
$\mathbb F$-vector space
$U_0 \otimes
\mathbb F \lbrack a^{\pm 1}, b^{\pm 1}, c^{\pm 1}\rbrack $:
\begin{gather*}
y^i \Lambda^j\otimes a^r b^s c^t, \qquad i,r,s,t\in \mathbb Z, \qquad j\in \mathbb N.
\end{gather*}
\end{Lemma}
\begin{proof}
By Lemma
\ref{lem:u0facts} and the construction.
\end{proof}

 Let $\lbrace \lambda_i\rbrace_{i=1}^5 $ denote mutually
commuting indeterminates.

\begin{Corollary}
\label{cor:polyalg}
There exists an $\mathbb F$-algebra isomorphism
\begin{gather*}
U_0 \otimes
\mathbb F \big\lbrack a^{\pm 1}, b^{\pm 1}, c^{\pm 1}\big\rbrack  \to
\mathbb F\big\lbrack \lambda_1,
\lambda_2^{\pm 1},
\lambda_3^{\pm 1},
\lambda_4^{\pm 1},
\lambda_5^{\pm 1} \big\rbrack
\end{gather*}
that sends
\begin{gather*}
\Lambda \otimes 1 \mapsto \lambda_1, \quad
y\otimes 1 \mapsto \lambda_2 ,\quad
1\otimes a \mapsto \lambda_3,
\quad
1\otimes b \mapsto \lambda_4,
\quad
1\otimes c \mapsto \lambda_5.
\end{gather*}
\end{Corollary}
\begin{proof}
Immediate from
Lemma
\ref{lem:u0tensor}.
\end{proof}

\begin{Definition}
\label{def:tildeproj}
\rm For $n \in \mathbb Z$  consider the
map
\begin{gather*}
\pi_n \otimes 1: \quad
\genfrac{}{}{0pt}{}{
U \otimes
\mathbb F \lbrack a^{\pm 1}, b^{\pm 1}, c^{\pm 1}\rbrack \; \to\;
U\otimes
\mathbb F \lbrack a^{\pm 1}, b^{\pm 1}, c^{\pm 1}\rbrack,
}
{\qquad  u \otimes f  \quad \qquad
\mapsto
\qquad \quad
\pi_n(u)\otimes f.}
\end{gather*}
The map
$\pi_n \otimes 1$
acts as the identity on the $n$-homogeneous component of
$U\otimes
\mathbb F \lbrack a^{\pm 1}, b^{\pm 1}, c^{\pm 1}\rbrack$
and  zero on all  other homogeneous components of
$U \otimes
\mathbb F \lbrack a^{\pm 1}, b^{\pm 1}, c^{\pm 1}\rbrack$.
Therefore
$\pi_n \otimes 1$ is the
projection from
$U\otimes
\mathbb F \lbrack a^{\pm 1}, b^{\pm 1}, c^{\pm 1}\rbrack$
onto its $n$-homogeneous component.
We abbreviate
$\tilde \pi_n =
\pi_n \otimes 1$.
So for $v \in
U \otimes
\mathbb F \lbrack a^{\pm 1}, b^{\pm 1}, c^{\pm 1}\rbrack$
the element
$\tilde \pi_n(v)$ is the homogeneous component of
$v$ that has degree $n$.
\end{Definition}

\begin{Lemma}
\label{lem:abctable}
In the table below we list some
elements $v \in
U \otimes
\mathbb F \lbrack a^{\pm 1}, b^{\pm 1}, c^{\pm 1}\rbrack $.
For each element $v$ we display the homogeneous
component $\tilde \pi_n(v)$ for $-1 \leq n \leq 1$.
All other homogeneous components
of $v$ are zero.
\begin{center}
\begin{tabular}{@{}c@{\,\,}|@{\,\,}c@{\,\,}c@{\,\,}c@{}}
$v$ &  $\tilde \pi_{-1}(v)$ & $\tilde \pi_0(v)$  & $\tilde \pi_1(v)$
\\
\hline
\hline
$A^\natural$ & $0$ & $y \otimes a^{-1} + y^{-1} \otimes a$ &
$\nu_z \otimes bc^{-1}-q^{-1} \nu_z y^{-1}  \otimes a$\tsep{2pt}
\\
$B^\natural $ & $\nu_x \otimes a^{-1}c-q^{-1} y^{-1} \nu_x \otimes b^{-1}$
& $y\otimes b + y^{-1} \otimes b^{-1}$\tsep{2pt}
& $0$
\\
$C^\natural$ &
$q^{-2}y^{-2} \nu_x  \otimes ab^{-1} $
&
  $y^{-1} \otimes (c+c^{-1})
+y^{-1}\Lambda \otimes ab^{-1}$
& $q^{-2}\nu_z y^{-2}\otimes ab^{-1}$\tsep{2pt}
\\
&
$-q^{-1} y^{-1} \nu_x \otimes c$
& $
 -(q+q^{-1})y^{-2}\otimes ab^{-1}$
&
$-q^{-1}\nu_z y^{-1} \otimes c^{-1}$
\end{tabular}
        \end{center}
Moreover each of
$\alpha^\natural$, $\beta^\natural$,
$\gamma^\natural$,
$\Omega^\natural $
is homogeneous with degree zero.
\end{Lemma}
\begin{proof}
The elements
$A^\natural$, $B^\natural$, $C^\natural$
are from 
(\ref{eq:aform})--(\ref{eq:cform})
and
$\alpha^\natural$, $\beta^\natural$, $\gamma^\natural$
are from 
(\ref{eq:amaincom1})--(\ref{eq:cmaincom3}).
Moreover $\Omega^\natural$ is from 
(\ref{eq:omimage}).
Evaluate these lines
using
Lemma~\ref{lem:firsttable}
and the fact that each of
$
1\otimes a^{\pm 1}$,
$1\otimes b^{\pm 1}$,
$1\otimes c^{\pm 1}$ is homogeneous of degree~0.
\end{proof}

The following def\/inition is for  notational convenience.

\begin{Definition}
\label{lem:commf}
Def\/ine $R$ and $L$ by
\begin{gather*}
R
 =
\tilde \pi_{1}(A^\natural)
=
\nu_z \otimes bc^{-1}-q^{-1} \nu_z y^{-1}\otimes a,
\\
L
=
\tilde \pi_{-1}(B^\natural)
=
\nu_x \otimes a^{-1}c-q^{-1} y^{-1} \nu_x \otimes b^{-1}.
\end{gather*}
Further def\/ine
\begin{gather}
\theta = y^{-1} \otimes a,
  \qquad
\vartheta = y^{-1} \otimes b^{-1}.
\label{eq:thetadef}
\end{gather}
\end{Definition}

\begin{Lemma}
\label{lem:comm}
The elements $R$, $L$, $\theta$, $\vartheta$ from
Definition~{\rm \ref{lem:commf}}
are all nonzero. Moreover
\begin{alignat}{3}
& R \vartheta = q^2 \vartheta R,
\qquad &&
L \theta = q^{-2} \theta L, &
\label{eq:rl1}
\\
&\tilde \pi_{0}(A^\natural) = \theta + \theta^{-1},
\qquad && \tilde \pi_{0}(B^\natural) = \vartheta + \vartheta^{-1}, &
\label{eq:rl2}
\\
&
\tilde \pi_{1}(C^\natural) =
-q^{-1} R \vartheta,
\qquad &&
\tilde \pi_{-1}(C^\natural) =
-q^{-1} \theta L.&
\label{eq:rl3}
\end{alignat}
\end{Lemma}
\begin{proof}
The f\/irst assertion
follows from
Lemma~\ref{lem:gradingbasis}.
(\ref{eq:rl1}) follows from~(\ref{eq:com2}) and
Def\/inition~\ref{lem:commf}.
(\ref{eq:rl2}),
(\ref{eq:rl3}) are readily checked using
the table in Lemma~\ref{lem:abctable}.
\end{proof}

We now give two lemmas of a slightly technical
nature.

\begin{Lemma}
\label{lem:aproj}
For an integer $i \geq 0$ the homogeneous components of
$
(A^\natural)^i$,
$(B^\natural)^i$,
$(C^\natural)^i$ are described as follows.

 $(A^\natural)^i$:
The homogeneous component of degree $n$ is zero unless
$0\leq n \leq i$.
The homogeneous component of degree $0$ is
$(\theta + \theta^{-1})^i$
and the homogeneous component of degree $i$ is~$R^i$.

 $(B^\natural)^i$:
The homogeneous component of degree $n$ is zero unless
$-i\leq n \leq 0$.
The homogeneous component of degree $-i$ is
$L^i$
and the
homogeneous component of degree $0$ is
$(\vartheta+ \vartheta^{-1})^i$.

$(C^\natural)^i$:
The homogeneous component of degree $n$ is zero unless
$-i\leq n \leq i$.
The homogeneous component of degree $-i$ is
$(-1)^iq^{i^2}L^i \theta^i$ and the
homogeneous component of degree $i$ is
$(-1)^iq^{-i^2}R^i \vartheta^i$.
\end{Lemma}
\begin{proof}
This is readily checked using
Lemma~\ref{lem:abctable} and
Lemma~\ref{lem:comm}.
\end{proof}

Using Lemma~\ref{lem:comm}
and Lemma~\ref{lem:aproj} we routinely obtain the following result.

\begin{Lemma}
\label{lem:threeproj}
Fix nonnegative integers $i$, $j$, $k$
and consider the homogeneous components of
$(A^\natural)^i(B^\natural)^j(C^\natural)^k$.
The homogeneous component of degree $n$ is zero unless
$-j-k\leq n \leq i+k$.
The homogeneous component of degree $-j-k$ is
\begin{gather*}
(-1)^k q^{k^2}
L^{j+k}
\big(
q^{2j+2k} \theta
+
q^{-2j-2k} \theta^{-1}
\big)^i
\theta^k.
\end{gather*}
The homogeneous component of degree $i+k$ is
\begin{gather*}
(-1)^k q^{-k^2}
R^{i+k}
\big(
q^{-2k}\vartheta +
q^{2k} \vartheta^{-1}
\big)^j
\vartheta^k.
\end{gather*}
\end{Lemma}

\section{Some results concerning algebraic independence}\label{section8}

In this section we establish some results about
algebraic independence that will be used in the proof of
Theorem~\ref{thm:main3}.

Let $\lbrace x_i\rbrace_{i=1}^4$ denote mutually
commuting indeterminates.
Motivated by the form of
(\ref{eq:omimage}) and
(\ref{eq:amaincom1})--(\ref{eq:cmaincom3}) we consider
the following
elements  in
$\mathbb F\lbrack x_1, x_2, x_3, x_4\rbrack$:
\begin{gather}
 y_1 =
x_1x_2x_3x_4
+x_1^2+x_2^2+x_3^2+x_4^2,
\label{eq:y1def}
\\
y_2 = x_1x_2+x_3x_4,
\qquad
y_3 = x_1x_3+x_2x_4,
\qquad
y_4 = x_1x_4+x_2x_3.
\label{eq:y234def}
\end{gather}

\begin{Lemma}
\label{lem:fiveindep}
The elements
 $\lbrace y_i\rbrace_{i=1}^4$ in
\eqref{eq:y1def},
\eqref{eq:y234def}
 are algebraically
independent over $\mathbb F$.
\end{Lemma}
\begin{proof}
The following is a basis for
the $\mathbb F$-vector space
$\mathbb F\lbrack x_1, x_2, x_3, x_4\rbrack$:
\begin{gather}
x_1^hx_2^ix_3^jx_4^k, \qquad h,i,j,k\in \mathbb N.
\label{eq:x1basis}
\end{gather}
An element
$x_1^hx_2^ix_3^jx_4^k$
in
the basis
(\ref{eq:x1basis}) will be called a {\it monomial}.
The {\it rank} of this monomial is def\/ined
to be $2h+i+j+k$.
For example,
consider the monomials that make up
$y_1$ in
(\ref{eq:y1def}).
The monomial
$x_1x_2x_3x_4$ has rank 5, and
the remaining monomials
$x_1^2$, $x_2^2$, $x_3^2$, $x_4^2$ have rank
$4$, $2$, $2$, $2$ respectively.
For $2 \leq i \leq 4$
consider the two monomials that make up $y_i$ in~(\ref{eq:y234def}). In each case
the monomial involving $x_1$ has
rank~3 and the other monomial
has rank~2.
To prove the lemma,
it suf\/f\/ices to show that the following elements
are linearly independent over $\mathbb F$:
\begin{gather}
y_1^ry_2^sy_3^ty_4^u, \qquad r,s,t,u\in \mathbb N.
\label{eq:yybasis}
\end{gather}
Given integers $r,s,t,u\geq 0$
 write
$y_1^ry_2^sy_3^ty_4^u$
as a linear combination of monomials:
\begin{gather*}
y_1^ry_2^sy_3^ty_4^u
=
(x_1x_2x_3x_4
+x_1^2+x_2^2+x_3^2+x_4^2)^r
(x_1x_2+x_3x_4)^s
(x_1x_3+x_2x_4)^t
(x_1x_4+x_2x_3)^u
\\
\phantom{y_1^ry_2^sy_3^ty_4^u}{} =
(x_1x_2x_3x_4)^r
(x_1x_2)^s
(x_1x_3)^t
(x_1x_4)^u
 + {\mbox{sum of monomials that have lower rank}}
\\
\phantom{y_1^ry_2^sy_3^ty_4^u}{} =
x_1^{r+s+t+u}
x_2^{r+s}
x_3^{r+t}
x_4^{r+u}
  +  {\mbox{sum of monomials that have lower rank}}.
\end{gather*}
Let us call
the monomial
$x_1^{r+s+t+u}
x_2^{r+s}
x_3^{r+t}
x_4^{r+u} $ the {\it leading monomial} for
$y_1^ry_2^sy_3^ty_4^u$. Given a~monomial $x_1^hx_2^ix_3^jx_4^k$ in the basis
(\ref{eq:x1basis}), consider the following
system of
linear equations in the unknowns $r$, $s$, $t$, $u$:
\begin{gather*}
r+s+t+u=h,
\qquad r+s=i,
\qquad r+t=j,
\qquad r+u=k.
\end{gather*}
Over the rational f\/ield $\mathbb Q$ this system has a unique
solution
\begin{gather*}
 r = \frac{i+j+k-h}{2},
\qquad
s = \frac{h+i-j-k}{2},
\qquad
 t = \frac{h-i+j-k}{2},
  \qquad
u = \frac{h-i-j+k}{2}.
\end{gather*}
Therefore
 $x_1^hx_2^ix_3^jx_4^k$
is the leading monomial
for at most one element
of (\ref{eq:yybasis}).
By these comments the elements
(\ref{eq:yybasis}) are linearly independent over $\mathbb F$.
The result follows.
\end{proof}

Recall the commutative algebra
$U_0\otimes
\mathbb F \lbrack a^{\pm 1}, b^{\pm 1}, c^{\pm 1}\rbrack $
from
Lemma \ref{lem:u0tensor}, and the element
$\theta = y^{-1}\otimes a$ from
(\ref{eq:thetadef}).

\begin{Proposition}
\label{lem:algind}
The following elements of
$U_0\otimes
\mathbb F \lbrack a^{\pm 1}, b^{\pm 1}, c^{\pm 1}\rbrack $
are algebraically independent over $\mathbb F$:
\begin{gather*}
\theta, \quad \Omega^\natural,
\quad
\alpha^\natural,
\quad
\beta^\natural,
\quad
\gamma^\natural.
\end{gather*}
\end{Proposition}
\begin{proof}
By
 Corollary
\ref{cor:polyalg} the following are algebraically independent over
$\mathbb F$:
\begin{gather*}
y \otimes 1, \quad
\Lambda \otimes 1,
\quad
1 \otimes a,
\quad
1 \otimes b,
\quad
1 \otimes c.
\end{gather*}
Therefore the following are
algebraically independent over
$\mathbb F$:
\begin{gather*}
y^{-1} \otimes a, \quad
\Lambda \otimes 1,
\quad
1 \otimes a,
\quad
1 \otimes b,
\quad
1 \otimes c.
\end{gather*}
Therefore the following are
algebraically independent over
$\mathbb F$:
\begin{gather}
y^{-1} \otimes a, \quad
\Lambda \otimes 1,
\quad
1 \otimes \big(a+a^{-1}\big),
\quad
1 \otimes \big(b+b^{-1}\big),
\quad
1 \otimes \big(c+c^{-1}\big).
\label{eq:seq5}
\end{gather}
Abbreviate the sequence
(\ref{eq:seq5})
by
$X_0$, $X_1$, $X_2$, $X_3$, $X_4$.
By
Lemma
\ref{lem:fiveindep}
the following are
algebraically independent over
$\mathbb F$:
\begin{gather*}
 X_0,
\qquad
X_1X_2X_3X_4
+X_1^2+X_2^2+X_3^2+X_4^2,
\\
 X_1X_2+X_3X_4,
\qquad
 X_1X_3+X_2X_4,
\qquad
 X_1X_4+X_2X_3.
\end{gather*}
The
above f\/ive elements are
\begin{gather*}
\theta,
 \quad
 1\otimes \big(q+q^{-1}\big)^2 -
 \Omega^\natural,
\quad
\alpha^\natural,
\quad
\beta^\natural,
\quad
\gamma^\natural,
\end{gather*}
respectively. The result follows.
\end{proof}

 Recall the element
$\vartheta = y^{-1} \otimes b^{-1}$ from~(\ref{eq:thetadef}).

\begin{Proposition}
\label{lem:algind2}
The following elements of
$U_0\otimes
\mathbb F \lbrack a^{\pm 1}, b^{\pm 1}, c^{\pm 1}\rbrack $
are algebraically independent over $\mathbb F$:
\begin{gather*}
\vartheta, \quad \Omega^\natural,
\quad
\alpha^\natural,
\quad
\beta^\natural,
\quad
\gamma^\natural.
\end{gather*}
\end{Proposition}
\begin{proof}
Similar to the proof of
Proposition
\ref{lem:algind}.
\end{proof}

\section{The proof of Theorem~\ref{thm:main3}}\label{section9}

In this section we prove
Theorem
\ref{thm:main3}.

\begin{Lemma}[\protect{\cite[Theorem 7.5]{uaw}}]
\label{lem:altbasisdelta}
The following is a basis for the $\mathbb F$-vector space
$\Delta$:
\begin{gather}
A^iB^jC^k\Omega^\ell \alpha^r \beta^s \gamma^t, \qquad
i,j,k,\ell, r,s,t\in \mathbb N,
\qquad   ijk=0.
\label{eq:truncbasis}
\end{gather}
\end{Lemma}

 We will be discussing
 the coef\/f\/icients
when an element of $\Delta$ is written as a linear
combination of the basis elements~(\ref{eq:truncbasis}). To facilitate this discussion
we def\/ine a bilinear form
$\langle \,,\,\rangle : \Delta \times \Delta \to \mathbb F$
such that $\langle u,v\rangle = \delta_{u,v}$
for all elements $u$, $v$ in the basis
(\ref{eq:truncbasis}). The bilinear form
$\langle\,,\,\rangle$ is symmetric, and
the basis
(\ref{eq:truncbasis}) is
orthonormal with respect to
$\langle \,,\, \rangle $.
For $v \in \Delta$,
\begin{gather*}
v = \sum \langle v, A^iB^jC^k\Omega^\ell\alpha^r\beta^s\gamma^t \rangle
 A^iB^jC^k\Omega^\ell\alpha^r\beta^s\gamma^t,
\label{eq:truncbilsum}
\end{gather*}
where the sum is over all elements
 $A^iB^jC^k\Omega^\ell\alpha^r\beta^s\gamma^t$ in the basis
(\ref{eq:truncbasis}). In this sum there are f\/initely
many nonzero summands.

\begin{proof}[Proof of Theorem \ref{thm:main3}.]
Let $J\subseteq \Delta$ denote the kernel
of $\natural$.
 We show that $J=0$. To do this we
assume $J\not=0$ and get a contradiction.
Fix $0 \not=v \in J$.
Let $S=S(v)$ denote the set of
7-tuples $(i,j,k,\ell,r,s,t)$ of nonnegative integers
such that $ijk=0$ and $\langle v, A^iB^jC^k\Omega^\ell \alpha^r
\beta^s \gamma^t\rangle \not=0$.
By construction
\begin{gather*}
v = \sum_{(i,j,k,\ell,r,s,t) \in S}
\langle v,A^iB^jC^k\Omega^\ell \alpha^r \beta^s \gamma^t \rangle
A^iB^jC^k\Omega^\ell \alpha^r \beta^s \gamma^t.
\end{gather*}
In this equation we apply $\natural$
to both sides and
 get
\begin{gather}
0 = \sum_{(i,j,k,\ell,r,s,t) \in S}
\langle v,A^iB^jC^k\Omega^\ell \alpha^r \beta^s \gamma^t \rangle
(A^\natural)^i (B^\natural)^j (C^\natural)^k (\Omega^\natural)^\ell
(\alpha^\natural)^r
(\beta^\natural)^s (\gamma^\natural)^t.
\label{eq:homimage}
\end{gather}
For an element
$(i,j,k,\ell,r,s,t) \in S$ def\/ine its {\it height}
to be $i+k$ and its {\it depth} to be
$j+k$.
For all integers $n\geq 0$ let
$S^+_n$ (resp.~$S^-_n$) denote the set of elements in $S$ that have
height $n$ (resp.\ depth $n$). By construction
$\lbrace S^+_n\rbrace_{n=0}^\infty$
(resp.\
$\lbrace S^-_n\rbrace_{n=0}^\infty$)
is a partition of $S$.
We assume $v \not=0$ so~$S$ is nonempty.
Therefore $\lbrace S^+_n\rbrace_{n=0}^\infty$
(resp.\
$\lbrace S^-_n\rbrace_{n=0}^\infty$)
are not all
empty.
By construction~$S$ has f\/inite cardinality,
so
f\/initely many of
$\lbrace S^+_n\rbrace_{n=0}^\infty$
(resp.\
$\lbrace S^-_n\rbrace_{n=0}^\infty$)
are nonempty.
Def\/ine
$N= \max  \lbrace n | S^+_n \not=\varnothing
\rbrace$
and
$M= \max \lbrace n | S^-_n \not=\varnothing
\rbrace$.
By construction $S^+_N$ and $S^-_M$ are nonempty.
We now split the argument into two cases.

 Case $N\leq M$:
Recall the projection map
$\tilde \pi_{-M}$ from
Def\/inition~\ref{def:tildeproj}.
Apply
$\tilde \pi_{-M}$
to
each
side of
(\ref{eq:homimage}).
Pick
$(i,j,k,\ell,r,s,t) \in S$
and consider the corresponding summand  in
(\ref{eq:homimage}).
The image of this summand under
$\tilde \pi_{-M}$ is computed using
 Lemma
\ref{lem:threeproj},
and found to be zero unless
$(i,j,k,\ell,r,s,t) \in S^-_M$.
The result of the  computation is that
\begin{gather*}
 0=L^M \!\!\!\!\!\sum_{(i,j,k,\ell,r,s,t) \in S^-_M}\!\!\!\!\!
\langle v, A^i B^jC^k\Omega^\ell \alpha^r \beta^s \gamma^t \rangle
({-}1)^k  q^{k^2}
\big(
q^{2M}
\theta
\!+\!
q^{-2M}\theta^{-1}
\big)^i
\theta^k
(\Omega^\natural)^\ell
(\alpha^\natural)^r
(\beta^\natural)^s (\gamma^\natural)^t,
\end{gather*}
where
$L$, $\theta$ are from
Def\/inition
\ref{lem:commf}.
By Lemma
\ref{lem:comm}
$L\not=0$.
We mentioned above Corollary
\ref{cor:nozerodiv} that
$U \otimes
\mathbb F \lbrack a^{\pm 1}, b^{\pm 1}, c^{\pm 1}\rbrack $
contains no zero divisors.
Therefore
\begin{gather}
 0=\sum_{(i,j,k,\ell,r,s,t) \in S^-_M}
\langle v, A^iB^jC^k\Omega^\ell \alpha^r \beta^s \gamma^t \rangle
(-1)^k  q^{k^2}
\big(
q^{2M}
\theta
+
q^{-2M}\theta^{-1}
\big)^i
\theta^k
\nonumber\\
\hphantom{0=\sum_{(i,j,k,\ell,r,s,t) \in S^-_M}}{}
\times
(\Omega^\natural)^\ell
(\alpha^\natural)^r
(\beta^\natural)^s (\gamma^\natural)^t.
\label{eq:aboveeq}
\end{gather}
Consider the above equation.
We noted earlier that
$S^-_M$ is nonempty.
By construction
the scalar
$\langle v, A^iB^jC^k\Omega^\ell \alpha^r \beta^s \gamma^t \rangle$
is nonzero for all
$(i,j,k,\ell,r,s,t) \in S^-_M$.
For
$(i,j,k,\ell,r,s,t) \in S^-_M$, at least one of
$i$, $j$, $k$ is zero since $ijk=0$. Moreover $j+k=M$ and
$i+k\leq M$.
For these constraints on~$i$,~$j$,~$k$ the possible solutions
for $(i,j,k)$ are
\begin{gather*}
(0,0,M),  (0,1,M-1),
\dots,  (0,M-1,1),  (0,M,0),
(1,M,0), \ldots, (M-1,M,0),  (M,M,0).
\end{gather*}
For the above values of $(i,j,k)$ the corresponding
values of
$(
q^{2M}
\theta
+
q^{-2M}\theta^{-1}
)^i
\theta^k $ are
\begin{gather*}
\theta^M, \theta^{M-1}, \ldots, \theta,  1,
q^{2M}
\theta
+
q^{-2M}\theta^{-1},
 \ldots,
\big(q^{2M}
\theta
+
q^{-2M}\theta^{-1}\big)^{M-1},
\big(q^{2M}
\theta
+
q^{-2M}\theta^{-1}\big)^M.
\end{gather*}
The above line
contains a sequence of
Laurent polynomials in $\theta$.
(\ref{eq:linebelow})
below contains a~sequence of
Laurent polynomials in $\theta$.
These two sequences are bases for the same vector space.
\begin{gather}
\label{eq:linebelow}
\theta^M,  \theta^{M-1},  \ldots,  \theta,   1,
 \theta^{-1},
 \ldots,
 \theta^{1-M},
 \theta^{-M}.
\end{gather}
With the above comments in mind,
 equation
(\ref{eq:aboveeq}) gives a nontrivial $\mathbb F$-linear dependency
among
\begin{gather*}
\theta^h
(\Omega^\natural)^\ell (\alpha^\natural)^r (\beta^\natural)^s
(\gamma^\natural)^t,
  \qquad \ell, r,s,t \in \mathbb N,\qquad
h \in \mathbb Z, \qquad
-M \leq h \leq M.
\end{gather*}
In the above line we multiply each term by $\theta^M$ and
obtain
a nontrivial $\mathbb F$-linear dependency among
\begin{gather*}
\theta^h
(\Omega^\natural)^\ell (\alpha^\natural)^r
(\beta^\natural)^s (\gamma^\natural)^t,
 \qquad h,\ell, r,s,t \in \mathbb N ,
\qquad  h \leq 2M.
\end{gather*}
The above linear dependency contradicts
Proposition~\ref{lem:algind},
for the present case
 $N\leq M$.

Case $M\leq N$:
The argument is similar to the previous case. However
the details are slightly dif\/ferent so we will show them.
Apply $\tilde \pi_{N}$ to
each
side of
(\ref{eq:homimage}).
Pick
$(i,j,k,\ell,r,s,t) \in S$
and consider the corresponding summand  in
(\ref{eq:homimage}).
The image of this summand under
$\tilde \pi_{N}$ is computed using
 Lemma~\ref{lem:threeproj},
and found to be zero unless
$(i,j,k,\ell,r,s,t) \in S^+_N$.
The result of the  computation is that
\begin{gather*}
 0=R^N \sum_{(i,j,k,\ell,r,s,t) \in S^+_N}
\langle v, A^i B^jC^k\Omega^\ell \alpha^r \beta^s \gamma^t \rangle
(-1)^k  q^{-k^2}
\big(
q^{-2k}
\vartheta
+
q^{2k}\vartheta^{-1}
\big)^j
\vartheta^k
\\
\hphantom{0=R^N \sum_{(i,j,k,\ell,r,s,t) \in S^+_N}}{}
\times \
(\Omega^\natural)^\ell
(\alpha^\natural)^r
(\beta^\natural)^s (\gamma^\natural)^t,
\end{gather*}
where
$R$, $\vartheta$ are from
Def\/inition~\ref{lem:commf}.
By Lemma~\ref{lem:comm}
$R\not=0$.
We have seen that
$U \otimes
\mathbb F \lbrack a^{\pm 1}, b^{\pm 1}, c^{\pm 1}\rbrack $
contains no zero divisors.
Therefore
\begin{gather}
 0=\sum_{(i,j,k,\ell,r,s,t) \in S^+_N}
\langle v, A^iB^jC^k\Omega^\ell \alpha^r \beta^s \gamma^t \rangle
(-1)^k  q^{-k^2}
\big(
q^{-2k}
\vartheta
+
q^{2k}\vartheta^{-1}
\big)^j
\vartheta^k
\nonumber\\
\hphantom{0=\sum_{(i,j,k,\ell,r,s,t) \in S^+_N}}{}
\times
(\Omega^\natural)^\ell
(\alpha^\natural)^r
(\beta^\natural)^s (\gamma^\natural)^t.
\label{eq:secondeq}
\end{gather}
Consider the above equation.
We noted earlier that
$S^+_N$ is nonempty.
By construction
the scalar
$\langle v, A^iB^jC^k\Omega^\ell \alpha^r \beta^s \gamma^t \rangle$
is nonzero for all
$(i,j,k,\ell,r,s,t) \in S^+_N$.
For
$(i,j,k,\ell,r,s,t) \in S^+_N$, at least one of
$i$, $j$, $k$ is zero since $ijk=0$. Moreover $i+k=N$ and
$j+k\leq N$.
For these constraints on~$i$,~$j$,~$k$ the possible solutions
for $(i,j,k)$ are
\begin{gather*}
(0,0,N),  (1,0,N-1),
\dots,  (N-1,0,1),  (N,0,0),
(N,1,0), \ldots, (N,N-1,0),  (N,N,0).
\end{gather*}
For the above values of $(i,j,k)$ the corresponding
values of
$(
q^{-2k}
\vartheta
+
q^{2k}\vartheta^{-1}
)^j
\vartheta^k $ are
\begin{gather*}
\vartheta^N, \vartheta^{N-1}, \ldots, \vartheta,  1,
\vartheta
+
\vartheta^{-1},
 \ldots,
\big(
\vartheta
+
\vartheta^{-1}\big)^{N-1},
\big(
\vartheta
+
\vartheta^{-1}\big)^N.
\end{gather*}
The above line
contains a sequence of
Laurent polynomials in $\vartheta$.
(\ref{eq2:linebelow})
below contains a sequence of
Laurent polynomials in $\vartheta$.
These two sequences are bases for the same vector space.
\begin{gather}
\label{eq2:linebelow}
\vartheta^N,  \vartheta^{N-1},  \ldots,  \vartheta,   1,
 \vartheta^{-1},
 \ldots,
 \vartheta^{1-N},
 \vartheta^{-N}.
\end{gather}
By these comments the  equation
(\ref{eq:secondeq})
gives a
 nontrivial $\mathbb F$-linear dependency
among
\begin{gather*}
\vartheta^h
(\Omega^\natural)^\ell (\alpha^\natural)^r (\beta^\natural)^s
(\gamma^\natural)^t,
 \qquad \ell, r,s,t \in \mathbb N, \qquad
h \in \mathbb Z, \qquad   -N \leq h \leq N.
\end{gather*}
In the above line we multiply each term by $\vartheta^N$ and
obtain
a nontrivial $\mathbb F$-linear dependency among
\begin{gather*}
\vartheta^h
(\Omega^\natural)^\ell (\alpha^\natural)^r
(\beta^\natural)^s (\gamma^\natural)^t
, \qquad h,\ell, r,s,t \in \mathbb N , \qquad
 h \leq 2N.
\end{gather*}
This linear dependency
contradicts
Proposition~\ref{lem:algind2},
for the present case
 $M\leq N$.

Both cases yield a contradiction
under the assumption
that $J\not=0$.
Therefore $J=0$ so $\natural$ is injective.
\end{proof}

\section[Comments on the ${\rm  {PSL}}_2(\mathbb Z)$ action]{Comments on the $\boldsymbol{{\rm  PSL}_2(\mathbb Z)}$ action}\label{section10}

Consider the injection
$\natural:\Delta \to
U \otimes
\mathbb F \lbrack a^{\pm 1}, b^{\pm 1}, c^{\pm 1}\rbrack $
 from
Theorem~\ref{thm:main}
and Theorem~\ref{thm:main3}.
Below Def\/inition~\ref{def:abc} we showed how
${\rm  {PSL}}_2(\mathbb Z)$ acts
 on~$\Delta$ as a group of automorphisms.
This
${\rm  {PSL}}_2(\mathbb Z)$ action
induces a
${\rm  {PSL}}_2(\mathbb Z)$ action on the image
$\Delta^\natural$.
It is reasonable to ask whether this action extends to
a
${\rm  {PSL}}_2(\mathbb Z)$ action on
$U\otimes
\mathbb F \lbrack a^{\pm 1}, b^{\pm 1}, c^{\pm 1}\rbrack $ as
a group of automorphisms.
This extension does not quite work;
let us
examine what goes wrong.

Recall the generators $\rho$, $\sigma$
of
${\rm  {PSL}}_2(\mathbb Z)$
from
below Def\/inition~\ref{def:abc}.
We will f\/irst consider the mathematics around~$\sigma $.

One can readily check using
Def\/inition~\ref{def:chev}
that there exists an automorphism of~$U$
that sends
\begin{gather*}
e \mapsto f,  \qquad f \mapsto e,  \qquad k^{\pm 1}\mapsto k^{\mp 1}.
\end{gather*}
More generally, for any nonzero $\xi \in \mathbb F$ there exists
an automorphism of~$U$ that sends
\begin{gather*}
e \mapsto \xi f,  \qquad f \mapsto \xi^{-1}e,
\qquad k^{\pm 1}\mapsto k^{\mp 1}.
\end{gather*}
The above automorphism
swaps
$U_n$ and $U_{-n}$
for all $n \in \mathbb Z$, where
$\lbrace U_n\rbrace_{n \in \mathbb Z}$ is the
$\mathbb Z$-grading
of~$U$ from
below Lemma~\ref{lem:efkbasis}.

With these  comments in mind we
now consider the algebra
$U \otimes
\mathbb F \lbrack a^{\pm 1}, b^{\pm 1}, c^{\pm 1}\rbrack $.

\begin{Lemma}
\label{lem:sigmaef}
There exists an automorphism
$\tilde \sigma $ of
$U \otimes
\mathbb F \lbrack a^{\pm 1}, b^{\pm 1}, c^{\pm 1}\rbrack $
that sends
\begin{alignat}{4}
& e \otimes 1 \mapsto f\otimes a^{-1}b^{-1}c,
\qquad  &&
f \otimes 1 \mapsto e\otimes abc^{-1},
\qquad &&
k\otimes 1 \mapsto k^{-1}\otimes 1,&
\label{eq:ss1}
\\
&
1 \otimes a \mapsto 1 \otimes b,
\qquad &&
1 \otimes b \mapsto 1 \otimes a,
\qquad &&
1 \otimes c \mapsto 1 \otimes c.&
\label{eq:ss2}
\end{alignat}
Moreover $\tilde \sigma^2=1$.
\end{Lemma}
\begin{proof}
There exists an $\mathbb F$-algebra homomorphism
$\tilde \sigma_1: U \to
U \otimes
\mathbb F \lbrack a^{\pm 1}, b^{\pm 1}, c^{\pm 1}\rbrack $
that sends
\begin{gather*}
e \mapsto f \otimes a^{-1}b^{-1}c,
\qquad
f \mapsto e \otimes abc^{-1},
\qquad
k \mapsto k^{-1} \otimes 1.
\end{gather*}
There exists an $\mathbb F$-algebra homomorphism
$\tilde \sigma_2:
\mathbb F \lbrack a^{\pm 1}, b^{\pm 1}, c^{\pm 1}\rbrack
\to
U \otimes
\mathbb F \lbrack a^{\pm 1}, b^{\pm 1}, c^{\pm 1}\rbrack $
that sends
\begin{gather*}
 a \mapsto 1 \otimes b,
\qquad
 b \mapsto 1 \otimes a,
\qquad
 c \mapsto 1 \otimes c.
\end{gather*}
Note that $\tilde \sigma_1(u)\tilde \sigma_2(f)=
\tilde \sigma_2(f) \tilde \sigma_1(u)$
for all $u \in
U$
and $f \in
\mathbb F \lbrack a^{\pm 1}, b^{\pm 1}, c^{\pm 1}\rbrack$.
By these comments
the map
\begin{gather*}
\tilde \sigma: \quad
\genfrac{}{}{0pt}{}{
U \otimes
\mathbb F \lbrack a^{\pm 1}, b^{\pm 1}, c^{\pm 1}\rbrack \; \to\;
U \otimes
\mathbb F \lbrack a^{\pm 1}, b^{\pm 1}, c^{\pm 1}\rbrack,
}
{\qquad  u \otimes f  \;\; \qquad
\mapsto
\qquad \quad
\tilde \sigma_1(u)\tilde \sigma_2(f)}
\end{gather*}
is an $\mathbb F$-algebra homomorphism
that satisf\/ies (\ref{eq:ss1}), (\ref{eq:ss2}).
One checks $\tilde \sigma^2=1$. Therefore $\tilde \sigma$ is
invertible and hence an automorphism.
\end{proof}

 Recall the $\mathbb Z$-grading
of
$U\otimes
\mathbb F \lbrack a^{\pm 1}, b^{\pm 1}, c^{\pm 1}\rbrack $ from
above Lemma
\ref{lem:gradingbasis}.

\begin{Lemma} The automorphism $\tilde \sigma$
of
$U \otimes
\mathbb F \lbrack a^{\pm 1}, b^{\pm 1}, c^{\pm 1}\rbrack $
has the following effect on the $\mathbb Z$-grading.
For $n\in \mathbb Z$, $\tilde \sigma $
swaps the
homogeneous components with degree
$n$, $-n$.
\end{Lemma}
\begin{proof}
In
(\ref{eq:ss1}), (\ref{eq:ss2}) we gave the
action of $\tilde \sigma$ on some homogeneous generators for
$U\otimes
\mathbb F \lbrack a^{\pm 1}, b^{\pm 1}, c^{\pm 1}\rbrack $.
For each generator
its image under $\tilde \sigma$
is homogeneous. Moreover
the generator and its image have opposite
degree.
The result follows.
\end{proof}

 We now consider the automorphism
$\tilde \sigma $ of
$U \otimes
\mathbb F \lbrack a^{\pm 1}, b^{\pm 1}, c^{\pm 1}\rbrack $
from the point of view of the equitable presentation.

\begin{Lemma}
\label{lem:autform}
In the table below we display some elements
$v$ of
$U \otimes
\mathbb F \lbrack a^{\pm 1}, b^{\pm 1}, c^{\pm 1}\rbrack $.
For each element $v$ we display the image $\tilde \sigma(v)$
under the map
$\tilde \sigma$ from
Lemma~{\rm \ref{lem:sigmaef}}.
\begin{center}
\begin{tabular}{c| c }
$v$ &  $\tilde \sigma(v)$
\\
\hline
\hline
$x\otimes 1$ &
$y\otimes 1 + \nu_x\otimes  a^{-1}b^{-1}c$
\\
$y\otimes 1$ &
         $y^{-1}\otimes 1$\tsep{2pt}
\\
$z\otimes 1$ &
$y\otimes 1 + \nu_z \otimes abc^{-1}$
\\
\hline
$\nu_x\otimes 1$ &
$-q^{-1}\nu_z y^{-1} \otimes abc^{-1}$\tsep{2pt}
\\
$\nu_y\otimes 1$ &
$-q \nu_z y \otimes abc^{-1} +
y \Lambda \otimes 1 -(q+q^{-1}) y^2 \otimes 1
- qy\nu_x \otimes a^{-1}b^{-1}c$
\\
$\nu_z\otimes 1$ &
$-q^{-1}y^{-1} \nu_x \otimes a^{-1}b^{-1}c$
\\
\hline
$\Lambda \otimes 1 $
&
$\Lambda \otimes 1 $
\end{tabular}
        \end{center}
\end{Lemma}

\begin{proof}
The images of
$y \otimes 1$,
$\nu_x\otimes 1$,
$\nu_z\otimes 1$ are obtained from
Lemma~\ref{lem:sigmaef}, using
Note~\ref{note:efview} and $y=k$.
The images of
$x\otimes 1$ and $z\otimes 1$ are now obtained
using~(\ref{eq:xzform}).
The image of $\Lambda \otimes 1$
is found using
Lemma~\ref{lem:sigmaef}
and~(\ref{eq:lambdaform}).
The image of
$\nu_y\otimes 1$ is found using
row $\nu_y$ in the table of
Lemma~\ref{lem:firsttable}.
\end{proof}

\begin{Lemma}
\label{lem:abgam}
The
map $\tilde \sigma $ from
Lemma~{\rm \ref{lem:sigmaef}}
sends
\begin{alignat}{4}
& \alpha^\natural \mapsto \beta^\natural,
\qquad &&
\beta^\natural \mapsto \alpha^\natural,
\qquad &&
\gamma^\natural \mapsto \gamma^\natural,&
\label{eq:abcnat}
\\
& A^\natural \mapsto B^\natural,
\qquad &&
B^\natural \mapsto A^\natural.&&&
\label{eq:ABnat}
\end{alignat}
\end{Lemma}
\begin{proof}
To verify
(\ref{eq:abcnat}),
in the equations
(\ref{eq:amaincom1})--(\ref{eq:cmaincom3})
apply $\tilde \sigma $ to each side,
and evaluate the result using
Lemma~\ref{lem:sigmaef}
and the fact that $\tilde \sigma$ f\/ixes
$\Lambda \otimes 1$.
To verify~(\ref{eq:ABnat})
we refer to rows
$A^\natural $ and $B^\natural$ of the table
in Lemma~\ref{lem:abctable}.
To each term in those rows,
apply $\tilde \sigma$ and
 evaluate the
result using  Lemma~\ref{lem:autform}.
\end{proof}

\begin{Proposition}
\label{lem:dc}
The following diagram commutes:
\[
\begin{CD}
\Delta  @>\natural>>
U \otimes
\mathbb F \big\lbrack a^{\pm 1}, b^{\pm 1}, c^{\pm 1}\big\rbrack
           \\
	  @V \sigma VV                     @VV \tilde \sigma V \\
\Delta	  @>>\natural>
U \otimes
\mathbb F \big\lbrack a^{\pm 1}, b^{\pm 1}, c^{\pm 1}\big\rbrack
		    \end{CD}
		    \]
\end{Proposition}
\begin{proof}
The $\mathbb F$-algebra $\Delta $ is generated by
$A$, $B$, $C$. By this and
(\ref{eq:u3}) the $\mathbb F$-algebra
$\Delta$ is generated by
$A$, $B$, $\gamma$.
Recall that~$\sigma $
swap $A$, $B$ and f\/ixes~$\gamma$.
By
Lemma~\ref{lem:abgam}
$\tilde \sigma $
swaps
 $A^\natural$, $B^\natural$ and f\/ixes $\gamma^\natural$.
The result follows.
\end{proof}

 Proposition~\ref{lem:dc}
shows that the action
of $\sigma$ on  $\Delta^\natural$ extends
to an automorphism $\tilde \sigma$ of $U \otimes
\mathbb F \lbrack a^{\pm 1}, b^{\pm 1}, c^{\pm 1}\rbrack $ that
has order~2.
So far so good.
 We now turn to the mathematics around~$\rho$.

Consider  the following variation on
$U$.

\begin{Definition}
\label{def:umod}
Def\/ine the $\mathbb F$-algebra
$U'$ by generators $X$, $Y$, $Z$ and
relations
\begin{gather*}
\frac{qXY-q^{-1}YX}{q-q^{-1}} = 1,
\qquad
\frac{qYZ-q^{-1}ZY}{q-q^{-1}} = 1,
\qquad
\frac{qZX-q^{-1}XZ}{q-q^{-1}} = 1.
\end{gather*}
\end{Definition}

 The above presentation
of $U'$ resembles
the equitable presentation of~ $U$, except that the generator~$Y^{-1}$
 is missing.

By construction
there exists an $\mathbb F$-algebra homomorphism
$\iota :
U'
\to
U$
that sends
\begin{gather*}
X \mapsto x,
\qquad
Y \mapsto y,
\qquad
Z \mapsto z.
\end{gather*}
We will need the fact that
$\iota$
is injective.
We will supply a proof shortly.

\begin{Lemma}
\label{lem:xyzbasis}
The following is a basis for the $\mathbb F$-vector space
$U$:
\begin{gather*}
x^h y^i z^j, \qquad h,j \in \mathbb N ,\qquad   i \in \mathbb Z.
\end{gather*}
\end{Lemma}
\begin{proof}
For all $n \in \mathbb N$ let
$V_n$ denote the subspace of
$U$ spanned by those
elements
$e^rk^sf^t$ from~(\ref{eq:kef}) that satisfy
 $r+t=n$.
By
Lemma~\ref{lem:efkbasis}
the
sum $U=
\sum\limits_{n=0}^\infty V_n$ is direct.
For all $h,j\in \mathbb N$ and $i \in \mathbb Z$
let us write $x^hy^iz^j$ in terms of~$e$,~$k$,~$f$.
By Proposition~\ref{prop:equit} and
Note~\ref{note:ident},
\begin{gather*}
x =
k^{-1} - e k^{-1} q^{-1}\big(q-q^{-1}\big),
\qquad   y=k,
\qquad
z = k^{-1}+ f\big(q-q^{-1}\big).
\end{gather*}
Using this together with
$ke=q^2ek$ and $kf=q^{-2}fk$, we f\/ind
$x^hy^iz^j  \in \sum\limits_{n=0}^{h+j} V_n $
and
\begin{gather*}
x^hy^iz^j  -
(-1)^h \big(q-q^{-1}\big)^{h+j}q^{-h^2} e^h k^{i-h} f^j
\in \sum_{n=0}^{h+j-1} V_n.
\end{gather*}
The result follows since
$\lbrace e^hk^{i-h}f^j\,|\, h,j\in \mathbb N, \; i\in \mathbb Z\rbrace$
is a basis for
$U$.
\end{proof}

\begin{Lemma}
\label{lem:reducedbasis}
The following is a basis for the $\mathbb F$-vector space
$U'$:
\begin{gather}
X^h Y^i Z^j, \qquad h,i,j \in \mathbb N.
\label{eq:base}
\end{gather}
\end{Lemma}
\begin{proof}
Using the relations in
Def\/inition~\ref{def:umod} we routinely f\/ind that
the elements
(\ref{eq:base}) span
$U'$.
The elements~(\ref{eq:base}) are linearly independent,
since
their images under
$\iota$ are linearly independent by
Lemma~\ref{lem:xyzbasis}.
The result follows.
\end{proof}

\begin{Lemma}
\label{lem:injection}
The above map $\iota :
U' \to
U$
is injective.
\end{Lemma}
\begin{proof}\sloppy
For the basis vectors~(\ref{eq:base}) their images under
$\iota$ are linearly independent by
 Lem\-ma~\ref{lem:xyzbasis}.
\end{proof}

Consider the subalgebra of
$U$ generated by $x$, $y$, $z$.
This subalgebra is the image of
$U'$ under~$\iota$. Invoking
 Lemma
\ref{lem:injection} we identify
this subalgebra with $U'$
 via
$\iota$.
The elements $\nu_x$, $\nu_y$, $\nu_z$
are contained in~$U'$
by
Def\/inition~\ref{def:nx}, and~$\Lambda \in U'$
by
Lemma~\ref{lem:sixforms}.
The algebra
$U' \otimes
\mathbb F \lbrack a^{\pm 1}, b^{\pm 1}, c^{\pm 1}\rbrack $
is the subalgebra  of
$U
 \otimes
\mathbb F \lbrack a^{\pm 1}, b^{\pm 1}, c^{\pm 1}\rbrack $
generated by
\begin{gather*}
x \otimes 1, \quad
y \otimes 1, \quad
z \otimes 1, \quad
1 \otimes a^{\pm 1}, \quad
1 \otimes b^{\pm 1}, \quad
1 \otimes c^{\pm 1}.
\end{gather*}
The next result clarif\/ies how
$U' \otimes
\mathbb F \lbrack a^{\pm 1}, b^{\pm 1}, c^{\pm 1}\rbrack $
is related to
$U \otimes
\mathbb F \lbrack a^{\pm 1}, b^{\pm 1}, c^{\pm 1}\rbrack $.

\begin{Lemma}
\label{lem:clarify}
The following is a basis for
the $\mathbb F$-vector space
$U\otimes
\mathbb F \lbrack a^{\pm 1}, b^{\pm 1}, c^{\pm 1}\rbrack $:
\begin{gather*}
x^hy^iz^j \otimes a^r b^s c^t, \qquad h,j\in \mathbb N, \qquad
i,r,s,t \in \mathbb Z.
\end{gather*}
The following is a basis for the
$\mathbb F$-vector space
$U'\otimes
\mathbb F \lbrack a^{\pm 1}, b^{\pm 1}, c^{\pm 1}\rbrack $:
\begin{gather*}
x^hy^iz^j \otimes a^r b^s c^t, \qquad h,i,j\in \mathbb N, \qquad
r,s,t \in \mathbb Z.
\end{gather*}
\end{Lemma}
\begin{proof}
Use Lemma
\ref{lem:xyzbasis} and
Lemma
\ref{lem:reducedbasis}.
\end{proof}

By (\ref{eq:aform})--(\ref{eq:cform})
the subalgebra
$U'\otimes
\mathbb F \lbrack a^{\pm 1}, b^{\pm 1}, c^{\pm 1}\rbrack $
contains $A^\natural$, $B^\natural$, $C^\natural$.
Therefore
$U'\otimes
\mathbb F \lbrack a^{\pm 1}, b^{\pm 1}, c^{\pm 1}\rbrack $
contains  $\Delta^\natural$.
 Consequently we may view
$\natural $ as an injection
$\natural: \Delta \to
U'\otimes
\mathbb F \lbrack a^{\pm 1}, b^{\pm 1}, c^{\pm 1}\rbrack $.

By Def\/inition
\ref{def:umod} there exists an automorphism of
$U'$
that sends
$(x,y,z)$ to
$(y,z,x)$.
This automorphism f\/ixes $\Lambda$ in view of
Lemma
\ref{lem:sixforms}.
There is also an automorphism of
$\mathbb F \lbrack a^{\pm 1}, b^{\pm 1}, c^{\pm 1}\rbrack $
that sends
 $(a,b,c)$ to $(b,c,a)$.
Combining these automorphisms we obtain the following.

\begin{Lemma}
\label{lem:rhodef}
There exists an automorphism $\tilde \rho$ of
$U' \otimes
\mathbb F \lbrack a^{\pm 1}, b^{\pm 1}, c^{\pm 1}\rbrack $
that sends
\begin{alignat*}{4}
&
x \otimes 1 \mapsto y\otimes 1,
\qquad &&
y \otimes 1 \mapsto z\otimes 1,
\qquad &&
z \otimes 1 \mapsto x\otimes 1,&
\\
&
1 \otimes a \mapsto 1\otimes b,
\qquad &&
1 \otimes b \mapsto 1\otimes c,
\qquad &&
1 \otimes c\mapsto 1\otimes a. &
\end{alignat*}
Moreover $\tilde \rho^3=1$.
\end{Lemma}

\begin{Lemma}
\label{lem:rhomove}
The automorphism $\tilde \rho$ from Lemma~{\rm \ref{lem:rhodef}} sends
\begin{gather*}
A^\natural \mapsto B^\natural,
\qquad
B^\natural \mapsto C^\natural,
\qquad
C^\natural \mapsto A^\natural.
\end{gather*}
\end{Lemma}
\begin{proof}
Evaluate
(\ref{eq:aform})--(\ref{eq:cform})
using
Lem\-ma~\ref{lem:rhodef}.
\end{proof}

\begin{Lemma}
The automorphism $\tilde \rho$ from Lemma~{\rm \ref{lem:rhodef}} fixes $\Lambda \otimes 1$ and sends
\begin{gather*}
\alpha^\natural \mapsto \beta^\natural,
\qquad
\beta^\natural \mapsto \gamma^\natural,
\qquad
\gamma^\natural \mapsto \alpha^\natural.
\end{gather*}
\end{Lemma}
\begin{proof}
This is routinely checked using
(\ref{eq:amaincom1})--(\ref{eq:cmaincom3})
and the comment about $\Lambda $ above
Lem\-ma~\ref{lem:rhodef}.
\end{proof}

\begin{Proposition}
\label{lem:dcrho}
The following diagram commutes:
\[
\begin{CD}
\Delta  @>\natural>>
U' \otimes
\mathbb F \lbrack a^{\pm 1}, b^{\pm 1}, c^{\pm 1}\rbrack
           \\
	  @V \rho VV                     @VV \tilde \rho V \\
\Delta	  @>>\natural>
U' \otimes
\mathbb F \lbrack a^{\pm 1}, b^{\pm 1}, c^{\pm 1}\rbrack
		    \end{CD}
		    \]
\end{Proposition}
\begin{proof}
The algebra $\Delta$ is generated by $A$, $B$, $C$.
Recall that $\rho$ cyclically
permutes $A$, $B$, $C$.
By Lemma~\ref{lem:rhomove}
$\tilde \rho$ of
cyclically
permutes $A^\natural$, $B^\natural$, $C^\natural$. The
result follows.
\end{proof}

Proposition~\ref{lem:dcrho}
shows that the action
of $ \rho$ on  $\Delta^\natural$ extends
to an automorphism $\tilde \rho$ of $U' \otimes
\mathbb F \lbrack a^{\pm 1}, b^{\pm 1}, c^{\pm 1}\rbrack $ that has order 3.
We now show that the
 automorphism $\tilde \rho$ of $U' \otimes
\mathbb F \lbrack a^{\pm 1}, b^{\pm 1}, c^{\pm 1}\rbrack $
does not extend to an automorphism of
  $U\otimes
\mathbb F \lbrack a^{\pm 1}, b^{\pm 1}, c^{\pm 1}\rbrack $.
For the moment assume that such
an extension exists. Since it is an automorphism,
it sends invertible elements to invertible elements.
The element $y \otimes 1$ is invertible
in
  $U \otimes
\mathbb F \lbrack a^{\pm 1}, b^{\pm 1}, c^{\pm 1}\rbrack $,
with inverse
$y^{-1}\otimes 1$. The element $z\otimes 1$ is not invertible
in
  $U\otimes
\mathbb F \lbrack a^{\pm 1}, b^{\pm 1}, c^{\pm 1}\rbrack $,
since $z$ is not invertible in
 $U$ by
\cite[Lemma~3.5]{equit}. This gives a contradiction
since $\tilde \rho$ sends $y \otimes 1$ to
$z \otimes 1$.
Therefore
the automorphism
 $\tilde \rho$ of $U' \otimes
\mathbb F \lbrack a^{\pm 1}, b^{\pm 1}, c^{\pm 1}\rbrack $
does not extend to an automorphism of
  $U \otimes
\mathbb F \lbrack a^{\pm 1}, b^{\pm 1}, c^{\pm 1}\rbrack $.

 In the above discussion we failed
 to obtain a
${\rm  {PSL}}_2(\mathbb Z)$ action on
  $U \otimes
\mathbb F \lbrack a^{\pm 1}, b^{\pm 1}, c^{\pm 1}\rbrack $.
Perhaps we should
search instead for a ${\rm  {PSL}}_2(\mathbb Z)$ action on
  $U'\otimes
\mathbb F \lbrack a^{\pm 1}, b^{\pm 1}, c^{\pm 1}\rbrack $.
Our extensions $\tilde \sigma$ and $\tilde \rho$
do not give such an action,
for the following reason.

\begin{Lemma}
The automorphism
$\tilde \sigma$ of
  $U\otimes
\mathbb F \lbrack a^{\pm 1}, b^{\pm 1}, c^{\pm 1}\rbrack $
does not leave
  $U' \otimes
\mathbb F \lbrack a^{\pm 1}, b^{\pm 1}, c^{\pm 1}\rbrack $
invariant.
\end{Lemma}
\begin{proof}
By Lemma
\ref{lem:autform}
the map $\tilde \sigma$ sends
$y\otimes 1$ to $y^{-1}\otimes 1$.
By Lemma~\ref{lem:clarify}
the subalgebra
$U'\otimes
\mathbb F \lbrack a^{\pm 1}, b^{\pm 1}, c^{\pm 1}\rbrack $
contains
$y\otimes 1$ but not
$y^{-1}\otimes 1$. The result follows.
\end{proof}

\begin{Problem}
\label{prob:final}
\rm
Find an $\mathbb F$-algebra $\mathcal A$ with the following
properties:
\begin{enumerate}\itemsep=0pt
\item[$(i)$]
There exists an injection of  $\mathbb F$-algebras
$\sharp:
  U \otimes
\mathbb F \lbrack a^{\pm 1}, b^{\pm 1}, c^{\pm 1}\rbrack \to \mathcal A$.
\item[$(ii)$]
The algebra $\mathcal A$ has an automorphism $\hat \sigma$ of order~2 that
makes the following diagram commute:
\[
\begin{CD}
U \otimes
\mathbb F \lbrack a^{\pm 1}, b^{\pm 1}, c^{\pm 1}\rbrack
  @>\sharp>> \mathcal A
	 \\
	  @V \tilde \sigma VV                     @VV \hat \sigma V \\
U\otimes
\mathbb F \lbrack a^{\pm 1}, b^{\pm 1}, c^{\pm 1}\rbrack
  @>>\sharp> \mathcal A
		    \end{CD}
		    \]
\item[$(iii)$]
The algebra $\mathcal A$ has an automorphism $\hat \rho$ of order 3 that
makes the following diagram commute:
\[
\begin{CD}
U'\otimes
\mathbb F \lbrack a^{\pm 1}, b^{\pm 1}, c^{\pm 1}\rbrack
  @>\sharp>> \mathcal A
	 \\
	  @V \tilde \rho VV                     @VV \hat \rho V \\
U' \otimes
\mathbb F \lbrack a^{\pm 1}, b^{\pm 1}, c^{\pm 1}\rbrack
  @>>\sharp> \mathcal A
		    \end{CD}
		    \]
\item[$(iv)$] There does not exist a proper subalgebra
of $\mathcal A$
that satisf\/ies $(i)$--$(iii)$ above.
\end{enumerate}
\end{Problem}

 The signif\/icance of the above problem
is summarized below.
\begin{Proposition}
\label{lem:final}
Let $\mathcal A$ denote
an $\mathbb F$-algebra that satisfies the four
conditions of Problem~{\rm \ref{prob:final}}. Then
${\rm  {PSL}}_2(\mathbb Z)$ acts on
$\mathcal A$ as a group of automorphisms
such that $\rho$ acts as $\hat \rho$ and
$\sigma$ acts as $\hat \sigma$. Moreover
the following diagram commutes for all $g \in
{\rm  {PSL}}_2(\mathbb Z)$:
\[
\begin{CD}
 \Delta
  @>\natural>>
U \otimes \mathbb F \lbrack a^{\pm 1}, b^{\pm 1}, c^{\pm 1}\rbrack
  @>\sharp>>
  \mathcal A
	 \\
	  @V g VV  &&
  @VVgV
\\
		   \Delta
  @>> \natural>
U \otimes \mathbb F \lbrack a^{\pm 1}, b^{\pm 1}, c^{\pm 1}\rbrack
  @>>\sharp>
   \mathcal A
		    \end{CD}
		    \]
\end{Proposition}
\begin{proof}
The
${\rm  {PSL}}_2(\mathbb Z)$ action on
$\mathcal A$ exists by the construction.
Concerning the diagram,
without loss we may assume that $g=\sigma$
or
 $g=\rho$.
For $g=\sigma$  the diagram commutes
by
Proposition~\ref{lem:dc}
and Problem~\ref{prob:final}$(ii)$.
For $g=\rho$  the diagram commutes
by
Proposition~\ref{lem:dcrho}
and Problem~\ref{prob:final}$(iii)$.
\end{proof}

\pdfbookmark[1]{References}{ref}
\LastPageEnding

\end{document}